\documentclass{article}

\usepackage{amssymb,amsmath,latexsym}

\newtheorem{thm}{Theorem}[section]
\newtheorem{pro}[thm]{Proposition}
\newtheorem{cor}[thm]{Corollary}
\newtheorem{lem}[thm]{Lemma}

\begin{document}


\title{\bf On the Cohomology of Modular Lie Algebras}

\author{J\"org Feldvoss\thanks{E-mail address: \tt jfeldvoss@jaguar1.usouthal.edu}\\
{\small Department of Mathematics and Statistics}\\{\small University of South Alabama}\\
{\small Mobile, AL 36688--0002, USA}}

\date{\small Dedicated to Robert L. Wilson and James Lepowsky on the occasion of their 
60$^\mathrm{th}$ birthdays}

\maketitle


\begin{abstract}

\noindent In this paper we establish a connection between the cohomology of a modular 
Lie algebra and its $p$-envelopes. We also compute the cohomology of Zassenhaus algebras 
and their minimal $p$-envelopes with coefficients in generalized baby Verma modules and 
in simple modules over fields of characteristic $p>2$. 
\medskip

\noindent {\it 2000 Mathematics Subject Classification.\/} Primary 17B50, 17B55, 17B56;
Secondary 17B10, 17B20, 17B70
\medskip

\noindent {\it Key words and phrases.\/} Cohomology, $p$-envelope, truncated (co)induced
module, Zassenhaus algebra, generalized baby Verma module, simple module, central 
extension
\end{abstract}


\section*{Introduction}


This paper is mainly a survey on the cohomology of modular Lie algebras but there are 
also several new results. It is intended to be sufficiently self-contained to serve as 
a first introduction to this topic which plays an important role in modular representation 
theory. Another goal is to advertise the use of $p$-envelopes and truncated (co)induced 
modules in the representation theory of Lie algebras over fields of prime characteristic 
$p$. As a main example throughout the paper we consider a series of rank one Lie algebras 
of Cartan type, the Zassenhaus algebras (which are simple in characteristic $p>2$). In 
the following we will describe the contents of the paper in more detail.

In the first section we give some background material on $p$-envelopes which will be useful 
for the rest of the paper. We refer the reader to \cite{Str1,Str2,SF} for more details and 
most of the proofs. The second section is devoted to establishing a connection between the 
cohomology of a modular Lie algebra and its $p$-envelopes which is very similar to a factorization 
theorem of Hochschild and Serre in characteristic zero. In particular, we show that the 
vanishing of the cohomology of a modular Lie algebra is equivalent to the vanishing of the 
cohomology of any of its $p$-envelopes. Parts of this section are contained in an unpublished
diploma thesis of Katrin Legler \cite{Leg}. In the third section we reformulate a cohomological 
vanishing theorem of Dzhumadil'daev in terms of the universal $p$-envelope, a result that 
rests mainly on the trivial action of the universal $p$-envelope $\hat{L}$ of $L$ on the 
cohomology of $L$. We also derive several consequences which will be useful later in the 
paper. In particular, we generalize a cohomological vanishing theorem for $\mathbb{Z}$-graded 
Lie algebras due to Chiu and Shen. In the fourth section we introduce truncated (co)induced 
modules which are important in the classification of simple modules for certain classes 
of modular Lie algebras. We recall the important fact that truncated induced modules and 
truncated coinduced modules coincide up to a twist. In particular, duals of truncated 
(co)induced modules are again truncated (co)induced modules, respectively. Moreover, we
give a short proof of Shapiro's lemma for truncated induced modules. Mil'ner used his 
version of truncated induced modules in order to classify the simple modules of the 
Zassenhaus algebras \cite{Mil}. In this paper we restrict our attention to certain 
truncated induced modules of the Zassenhaus algebras, the so-called generalized baby 
Verma modules which are sufficient for computing the cohomology of Zassenhaus algebras 
with coefficients in simple modules. In the final section we reprove Dzhumadil'daev's 
cohomological reduction theorem for Zassenhaus algebras by using Shapiro's lemma for 
truncated induced modules. Then we apply this to compute the $1$-cohomology of Zassenhaus 
algebras with coefficients in generalized baby Verma modules and in simple modules over 
fields of characteristic $p>2$ without using Mil'ner's classification of the simple modules 
for the Zassenhaus algebras. As a consequence we determine the central extensions of the 
Zassenhaus algebras. We also apply the main result of Section 2 in order to compute the 
$1$-cohomology of the minimal $p$-envelopes of the Zassenhaus algebras with coefficients 
in generalized baby Verma modules and in simple modules as well as their central extensions. 
The results for the minimal $p$-envelopes of the Zassenhaus algebras seem to be new.

For the notation and fundamental results from the representation theory of modular Lie
algebras we refer the reader to \cite{SF} and \cite{Str2}.
\bigskip

\noindent {\it Acknowledgments.\/} The author would like to thank the referee for several 
useful comments. While preparing the final version of this paper Walter Michaelis made many 
valuable suggestions concerning the clarity and the style of the paper for which the author 
is very grateful.


\section{$p$-Envelopes}


Let $L$ be a Lie algebra over a field of prime characteristic $p$. A triple $(\mathfrak{L},
(\cdot)^{[p]},\iota)$ is called a {\it $p$-envelope\/} of $L$ if $\mathfrak{L}$ is a 
restricted Lie algebra with $p$-mapping $(\cdot)^{[p]}$ and $\iota$ is a Lie algebra 
monomorphism from $L$ into $\mathfrak{L}$ such that the restricted subalgebra $\langle
\iota(L)\rangle_p$ of $\mathfrak{L}$ generated by $\iota(L)$ coincides with $\mathfrak{L}$.

\begin{pro}\label{comm}
Let $L$ be a Lie algebra over a field of prime characteristic $p$ and let $(\mathfrak{L},
(\cdot)^{[p]},\iota)$ be a $p$-envelope of $L$. Then the following statements hold:
\begin{enumerate}
\item[{\rm (1)}] $[\mathfrak{L},\mathfrak{L}]\subseteq\iota(L)$.
\item[{\rm (2)}] If $\mathfrak{I}$ is a subspace of $\mathfrak{L}$ such that $\iota(L)
                 \subseteq\mathfrak{I}$, then $\mathfrak{I}$ is an ideal of $\mathfrak{L}$.
\end{enumerate}
\end{pro}

\noindent {\it Proof.\/} (1): Since $\langle\iota(L)\rangle_p=\mathfrak{L}$, it follows 
from \cite[Proposition 2.1.3(2)]{SF} that $$[\mathfrak{L},\mathfrak{L}]=[\langle\iota(L)
\rangle_p,\langle\iota(L)\rangle_p]=[\iota(L),\iota(L)]\subseteq\iota(L)\,.$$ 

(2): By using (1), we obtain that $[\mathfrak{L},\mathfrak{I}]\subseteq[\mathfrak{L},
\mathfrak{L}]\subseteq\iota(L)\subseteq\mathfrak{I}$.\quad $\Box$
\bigskip

\noindent Moreover, the following structural features are preserved by $p$-envelopes. In 
particular, $p$-envelopes of solvable (nilpotent, abelian) Lie algebras are solvable 
(nilpotent, abelian, respectively).

\begin{pro}\label{solv}
Let $L$ be a Lie algebra over a field of prime characteristic $p$ and let 
$(\mathfrak{L},(\cdot)^{[p]},\iota)$ be a $p$-envelope of $L$. Then the following 
statements hold:
\begin{enumerate}
\item[{\rm (1)}] $L$ is abelian if and only if $\mathfrak{L}$ is abelian.
\item[{\rm (2)}] $L$ is nilpotent if and only if $\mathfrak{L}$ is nilpotent.
\item[{\rm (3)}] $L$ is solvable if and only if $\mathfrak{L}$ is solvable.
\end{enumerate}
\end{pro}

\noindent {\it Proof.\/} By definition $\langle\iota(L)\rangle_p=\mathfrak{L}$ and
therefore the assertions are immediate consequences of \cite[Proposition 2.1.3(2)]{SF}.
\quad $\Box$
\bigskip

\noindent {\it Remark.\/} Note that every abelian Lie algebra is restrictable (e.g. by 
the trivial $p$-mapping) but neither nilpotent Lie algebras nor solvable Lie algebras 
are necessarily restrictable (cf.~\cite[Chapter I, \S4, Exercise 24]{Bou} as well as 
\cite[Example 3 on p.~72 and Exercise 2.2.9]{SF}).
\bigskip

A $p$-envelope $(\mathfrak{L}_{\min},(\cdot)^{[p]},\iota)$ of a finite-dimensional 
modular Lie algebra $L$ is called {\it minimal\/} if its dimension is minimal among
the dimensions of all $p$-envelopes of $L$. Any two minimal $p$-envelopes are 
isomorphic as ordinary Lie algebras (cf.~\cite[Theorem 2.5.8(1)]{SF}). Of importance
is the well-known fact that a $p$-envelope $(\mathfrak{L},(\cdot)^{[p]},\iota)$ of $L$ 
is minimal if and only if the center $C(\mathfrak{L})$ of $\mathfrak{L}$ is contained 
in $\iota(L)$ (cf.~\cite[Theorem 2.5.8(3)]{SF}).

It turns out that minimal $p$-envelopes of simple Lie algebras are simple as restricted
Lie algebras (but in general not simple as ordinary Lie algebras since a non-restrictable 
Lie algebra is isomorphic to a non-zero proper ideal of any of its $p$-envelopes) and 
minimal $p$-envelopes of semisimple Lie algebras are semisimple.

\begin{pro}\label{semsim}
Let $L$ be a finite-dimensional Lie algebra over a field of prime characteristic $p$ 
and let $(\mathfrak{L}_{\min},(\cdot)^{[p]},\iota)$ be a minimal $p$-envelope of $L$. 
Then the following statements hold:
\begin{enumerate}
\item[{\rm (1)}] If $L$ is simple, then $\mathfrak{L}_{\min}$ has no non-zero proper
                 $p$-ideals.
\item[{\rm (2)}] If $L$ is semisimple, then $\mathfrak{L}_{\min}$ is semisimple.
\end{enumerate}
\end{pro}

\noindent {\it Proof.\/} (1): Let $\mathfrak{I}$ be a $p$-ideal of $\mathfrak{L}_{\min}$.
Then $\mathfrak{I}\cap\iota(L)$ is an ideal of $\iota(L)$. But because $\iota(L)$ is simple,
either $\mathfrak{I}\cap\iota(L)=0$ or $\mathfrak{I}\cap\iota(L)=\iota(L)$. 

In the first case it follows from Proposition \ref{comm}(1) that $[\mathfrak{I},
\mathfrak{L}_{\min}]\subseteq\mathfrak{I}\cap[\mathfrak{L}_{\min},\mathfrak{L}_{\min}]
\subseteq\mathfrak{I}\cap\iota(L)=0$, i.e., that $\mathfrak{I}\subseteq C(\mathfrak{L}_{\min})
\subseteq\iota(L)$ which implies that $\mathfrak{I}=\mathfrak{I}\cap\iota(L)=0$.

In the second case one has $\iota(L)\subseteq\mathfrak{I}$. Since $\mathfrak{I}$ is a 
$p$-ideal of $\mathfrak{L}_{\min}$, one obtains that $\mathfrak{L}_{\min}=\langle\iota(L)
\rangle_p\subseteq\mathfrak{I}$, i.e., that $\mathfrak{I}=\mathfrak{L}_{\min}$.
\medskip

(2): Let $\mathfrak{I}$ be an abelian ideal of $\mathfrak{L}_{\min}$. Then $\mathfrak{I}
\cap\iota(L)$ is an abelian ideal of $\iota(L)$. But because $\iota(L)$ is semisimple,  
$\mathfrak{I}\cap\iota(L)=0$. 

As in the proof of (1), one has $\mathfrak{I}\subseteq C(\mathfrak{L}_{\min})\subseteq
\iota(L)$. Since $\iota(L)$ is semisimple and $\mathfrak{I}$ is an abelian ideal of 
$\iota(L)$, it readily follows that $\mathfrak{I}=0$.\quad $\Box$
\bigskip

We conclude this section by introducing the Zassenhaus algebras and their minimal 
$p$-envelopes. Let $\mathbb{F}$ be a field of prime characteristic $p$, let $m$ be 
a positive integer, and let $A(m)$ denote the subalgebra of the algebra of divided 
powers over $\mathbb{F}$ generated by $\{x^{(a)}\mid 0\le a\le p^m-1\}$ where 
$x^{(a)}x^{(b)}={a+b\choose a}x^{(a+b)}$ for non-negative integers $a$ and $b$. A 
derivation $d$ of $A(m)$ is called {\it special\/} if $d(x^{(0)})=0$ and $d(x^{(a)})
=x^{(a-1)}d(x^{(1)})$ for every $1\le a\le p^m-1$. The Lie algebra $W(m):=\{d\in
\mathrm{Der}_\mathbb{F}(A(m))\mid d\mbox{ is special}\}$ of special derivations of 
$A(m)$ for some positive integer $m$ is called a {\it Zassenhaus algebra\/}. It is 
well-known that $W(m)$ is simple if and only if $p>2$ (cf.~\cite[Theorem 4.2.4(1)]{SF}). 
Furthermore, $W(m)$ is restricted if and only if $m=1$ (cf.~\cite[Theorem 4.2.4(2)]{SF}). 
The {\it Witt algebra\/} $W(1)$ was the first non-classical simple Lie algebra 
discovered by Ernst Witt in the late 1930's. 

Let $\partial$ denote the derivative defined by
\begin{eqnarray*}
\partial(x^{(a)})=
\left\{
\begin{array}{cl}
0 & \mbox{if }a=0\\ 
x^{(a-1)} & \mbox{if }1\le a\le p^m-1\,. 
\end{array}
\right.
\end{eqnarray*}
Note that $W(m)$ is a free $A(m)$-module with basis $\partial$ (cf.~\cite[Proposition 
4.2.2(1)]{SF}). Set $e_i:=x^{(i+1)}\partial$ for any $-1\le i\le p^m-2$. Then 
\begin{eqnarray*}
[e_i,e_j]=
\left\{
\begin{array}{cl}
\left({i+j+1\choose j}-{i+j+1\choose i}\right)e_{i+j} & \mbox{if }-1\le i+j\le p^m-2\\
0 & \mbox{otherwise}
\end{array}
\right.
\end{eqnarray*}
and $$W(m)=\bigoplus_{i=-1}^{p^m-2}\mathbb{F}e_i$$ is a $\mathbb{Z}$-graded Lie algebra 
(cf.~\cite[Proposition 4.2.2(2) and (3)]{SF}). 

The subalgebra $$\mathfrak{b}(m):=\bigoplus_{i=0}^{p^m-2}\mathbb{F}e_i$$ is supersolvable
(i.e., every composition factor of the adjoint module of $\mathfrak{b}(m)$ is one-dimensional). 
Moreover, $\mathfrak{b}(m)$ is a restricted Lie algebra via the $p$-mapping defined by
\begin{eqnarray*}
e_i^{[p]}=
\left\{
\begin{array}{cl}
e_i & \mbox{if }i=0\\ 
(p-1)!e_{pi} & \mbox{if }i=p^t-1\mbox{ for some }1\le t\le m-1\\
0 & \mbox{otherwise} 
\end{array}
\right.
\end{eqnarray*}
(cf.~\cite[(3.2.2)]{Shu}). In particular, $\mathfrak{t}:=\mathbb{F}e_0$ is a one-dimensional 
torus of $\mathfrak{b}(m)$. Since $e_i^{[p]^2}=0$ for every $1\le i\le p^m-2$, the $p$-subalgebra
$$\mathfrak{u}(m):=\bigoplus_{i=1}^{p^m-2}\mathbb{F}e_i$$ is $p$-unipotent so that $\mathfrak{b}(m)
=\mathfrak{t}\oplus\mathfrak{u}(m)$ is strongly solvable as are the {\it Borel subalgebras\/} of 
Lie algebras of reductive groups.

Suppose that $L$ is a modular Lie algebra with $C(L)=0$. Let $\mathfrak{L}_{\min}$ denote 
the $p$-subalgebra of $\mathfrak{gl}(L)$ generated by $\mathrm{ad}(L)$. Then $(\mathfrak{L}_{\min},
(\cdot)^p,\mathrm{ad})$ is a minimal $p$-envelope of $L$ (cf.~\cite[p.~97]{SF}). 

Let us apply this to the Zassenhaus algebras $W(m)$. It is clear from the $\mathbb{Z}$-grading 
of $W(m)$ that $(\mathrm{ad}\,e_{-1})^{p^m}=0$. Hence $$\mathfrak{W}(m)=\bigoplus_{r=1}^{m-1}
\mathbb{F}e_{-1}^{[p]^r}\oplus W(m)$$ is a minimal $p$-envelope of $W(m)$. Here $e_i$ 
is identified with $\mathrm{ad}\,e_i$ for every $-1\le i\le p^m-2$ while $e_{-1}^{[p]^r}$ is 
identified with $(\mathrm{ad}\,e_{-1})^{p^r}$ for every $1\le r\le m-1$. We propose to call 
the restricted Lie algebra $\mathfrak{W}(m)$ a {\it restricted Zassenhaus algebra\/}. It 
follows from \cite[Theorem 4.2.4(1)]{SF} and Proposition \ref{semsim} that $\mathfrak{W}(m)$ 
is semisimple and has no non-zero proper $p$-ideals if $p>2$. In particular, the above 
$p$-mapping of $\mathfrak{W}(m)$ is the only one possible (cf.~\cite[Corollary 2.2.2(1)]{SF}). 
Note also that $\mathfrak{W}(m)$ is the semidirect product of a {\it cyclic\/} $p$-unipotent 
subalgebra with the ideal $W(m)$.


\section{Cohomology of $p$-Envelopes}


In this section we investigate the relationship between the ordinary cohomology
of finite-dimensional Lie algebras and their finite-dimensional $p$-envelopes.
The first lemma we need is a consequence of a long exact sequence due to Dixmier
\cite[Proposition 1]{Dix} (cf.~\cite[Korollar 3.7]{Leg}):

\begin{lem}\label{codim1}
Let $L$ be a finite-dimensional Lie algebra over a field $\mathbb{F}$ of arbitrary 
characteristic, let $M$ be a finite-dimensional $L$-module, and let $I$ be an ideal 
in $L$ of codimension $1$. Then, for every non-negative integer $n$, there is an 
isomorphism $$H^n(L,M)\cong H^n(I,M)^L\oplus H^{n-1}(I,M)^L$$ of $\mathbb{F}$-vector 
spaces. 
\end{lem}

\noindent {\it Proof.\/} Because $\dim_\mathbb{F}L/I=1$, there exists an element $x
\in L$ such that $L=I\oplus\mathbb{F}x$. Let $n$ be a non-negative integer. According to 
\cite[Proposition 1]{Dix}, the sequence \[H^{n-1}(I,M)\stackrel{u_{n-1}}{\longrightarrow}
H^{n-1}(I,M)\stackrel{s_n}{\longrightarrow} H^n(L,M)\stackrel{r_n}{\longrightarrow}H^n(I,M)
\stackrel{u_n}{\longrightarrow}H^n(I,M)\] is exact, where $u_n$ is induced by $\theta^n(x)
:C^n(I,M)\to C^n(I,M)$ (cf.~\cite[Chapter I, \S3, Exercise 12(a)]{Bou}), $s_n$ is induced 
by the connecting homomorphism, and $r_n$ is induced by the restriction mapping from 
$C^n(L,M)$ to $C^n(I,M)$. Hence $\mathrm{Ker}(s_n)=\mathrm{Im}(u_{n-1})=LH^{n-1}(I,M)$ 
and $\mathrm{Im}(r_n)=\mathrm{Ker}(u_n)=H^n(I,M)^L$. Consequently, one has the short 
exact sequence \[(*)\qquad 0\longrightarrow H^{n-1}(I,M)/LH^{n-1}(I,M)\longrightarrow 
H^n(L,M)\stackrel{r_n}{\longrightarrow}H^n(I,M)^L\longrightarrow 0\,.\] It follows from 
the first isomorphism theorem that $H^{n-1}(I,M)/\mathrm{Ker}(u_{n-1})\cong\mathrm{Im}
(u_{n-1})$ and thus $$H^{n-1}(I,M)\cong\mathrm{Ker}(u_{n-1})\oplus\mathrm{Im}(u_{n-1})
\cong H^{n-1}(I,M)^L\oplus LH^{n-1}(I,M)\,,$$ i.e., $H^{n-1}(I,M)/LH^{n-1}(I,M)\cong 
H^{n-1}(I,M)^L$. This in conjunction with the split short exact sequence $(*)$ yields 
the assertion.\quad $\Box$
\bigskip

\noindent {\it Remark.\/} Lemma \ref{codim1} can also be obtained from the Hochschild-Serre 
spectral sequence (see \cite[pp.~30--33]{Leg}). Our direct proof has the advantage that 
the mappings in the short exact sequence $(*)$ are known explicitly. 
\bigskip

Let $(\hat{L},(\cdot)^p,\hat{\iota})$ denote the {\it universal $p$-envelope\/} of a Lie 
algebra $L$ over a field of prime characteristic $p$ (cf.~\cite[Definition 2, p.~92]{SF}) 
and let $M$ be an $L$-module with corresponding representation $\rho:L\to\mathfrak{gl}(M)$. 
By virtue of the universal property of $\hat{L}$, there exists a unique restricted Lie 
algebra homomorphism $\hat{\rho}:\hat{L}\to\mathfrak{gl}(M)$ such that $\hat{\rho}\circ
\hat{\iota}=\rho$.

Now consider an arbitrary $p$-envelope $(\mathfrak{L},(\cdot)^{[p]},\iota)$ of $L$. In view 
of \cite[Proposition 2.5.6]{SF}, there exists a Lie algebra homomorphism $\phi:\mathfrak{L}
\to\hat{L}$ such that $\phi\circ\iota=\hat{\iota}$. Then $\hat{\rho}\circ\phi:\mathfrak{L}
\to\mathfrak{gl}(M)$ is a Lie algebra homomorphism such that $(\hat{\rho}\circ\phi)\circ
\iota=\rho$. In this way every $L$-module can be considered as an $\mathfrak{L}$-module
(cf.~\cite[Theorem 5.1.1]{SF}). Note that this representation of $\mathfrak{L}$ on $M$
is {\it not\/} unique but depends on the choice of the Lie algebra homomorphism $\phi:
\mathfrak{L}\to\hat{L}$.

It will be crucial for the proof of the main result of this section that the universal 
$p$-envelope of a modular Lie algebra $L$ acts trivially on the cohomology of $L$
(cf.~\cite[Lemma 2.10]{Leg}). Slightly more generally, we will need the following 
result: 

\begin{lem}\label{triv}
Let $L$ be a Lie algebra over a field of prime characteristic $p$ and let $M$ be an 
$L$-module. If $\iota(L)\subseteq\mathfrak{H}\subseteq\mathfrak{K}$ are subalgebras 
of a $p$-envelope $\mathfrak{L}$ of $L$, then $H^n(\mathfrak{H},M)$ is a trivial 
$\mathfrak{K}$-module for every non-negative integer $n$.
\end{lem}

\noindent {\it Proof.\/} Let $n$ be a non-negative integer. It is a consequence of 
\cite[Chapter I, \S3, Exercise 12(b), (3)]{Bou} that $H^n(L,M)$ is a trivial $L$-module 
for every non-negative integer $n$. Since $M$ is a restricted $\hat{L}$-module and 
the adjoint representation is restricted, one obtains that $H^n(L,M)$ is a restricted 
$\hat{L}$-module (cf.~\cite[Chapter I, \S3, Exercise 12(a)]{Bou} and \cite[Theorem 
5.2.7]{SF}). It follows that $\mathrm{Ann}_{\hat{L}}(H^n(L,M))$ is a $p$-subalgebra 
of $\hat{L}$ that contains $\hat{\iota}(L)$. Hence $\hat{L}=\langle\hat{\iota}(L)
\rangle_p=\mathrm{Ann}_{\hat{L}}(H^n(L,M))$, i.e., $H^n(L,M)$ is a trivial $\hat{L}$-module.

Now let $(\mathfrak{L},(\cdot)^{[p]},\iota)$ be an arbitrary $p$-envelope of $L$. Since 
the pullback functor from $\hat{L}$-modules to $\mathfrak{L}$-modules induced by $\phi$ 
commutes with tensor products and taking duals, the $\mathfrak{L}$-action on $H^n(L,M)$
factors through that of $\hat{L}$ and therefore is also trivial.

Finally, because of $\langle\mathfrak{H}\rangle_p=\langle\iota(L)\rangle_p=\mathfrak{L}$,
one can assume without loss of generality that $\mathfrak{H}=L$, hence that $\mathfrak{K}
\subseteq\mathfrak{L}$ acts trivially on $H^n(\mathfrak{H},M)$.\quad $\Box$
\bigskip

Now we derive dimension formulas for the cohomology of a finite-dimensional modular Lie 
algebra in terms of the cohomology of any of its finite-dimensional $p$-envelopes and 
vice versa (see \cite[Satz 3.9]{Leg}).

\begin{thm}\label{cohpenv}
Let $L$ be a finite-dimensional Lie algebra over a field $\mathbb{F}$ of prime characteristic 
$p$, let $M$ be a finite-dimensional $L$-module, and let $\mathfrak{L}$ be a $p$-envelope of 
$L$ with $k:=\dim_\mathbb{F}\mathfrak{L}/L<\infty$. Then the following statements hold:
\begin{enumerate}
\item[{\rm (1)}] $\dim_\mathbb{F}H^n(\mathfrak{L},M)=\sum\limits_{j=0}^k{k\choose j}
                 \dim_\mathbb{F}H^{n-j}(L,M)$ for every non-negative integer $n$.
\item[{\rm (2)}] $\dim_\mathbb{F}H^n(L,M)=\sum\limits_{\nu=0}^n(-1)^\nu{\nu+k-1\choose k-1}
                 \dim_\mathbb{F}H^{n-\nu}(\mathfrak{L},M)$ for every non-negative integer $n$.
\end{enumerate}
\end{thm}

\noindent {\it Proof.\/} (1): Choose a chain of subspaces $$\mathfrak{L}=\mathfrak{I}_k
\supset\mathfrak{I}_{k-1}\supset\cdots\supset\mathfrak{I}_1\supset\mathfrak{I}_0=L$$ with 
$\dim_\mathbb{F}\mathfrak{I}_l/\mathfrak{I}_{l-1}=1$ for every $1\le l\le k$. According 
to Proposition \ref{comm}(2), $\mathfrak{I}_l$ is an ideal of $\mathfrak{L}$ for every 
$0\le l\le k$. 

Let $n$ be a non-negative integer. It follows from Lemma \ref{codim1} and Lemma \ref{triv}
that $$(**)\qquad H^n(\mathfrak{I}_l,M)\cong H^n(\mathfrak{I}_{l-1},M)\oplus H^{n-1}
(\mathfrak{I}_{l-1},M)$$ for every $0\le l\le k$. 

Set $a_{n,l}:=\dim_\mathbb{F}H^n(\mathfrak{I}_l,M)$ for every non-negative integer $n$ and
every integer $0\le l\le k$. Then (1) can be written as $$a_{n,k}=\sum_{l=0}^k{k\choose l}
a_{n-l,0}\qquad\forall~n\in\mathbb{N}_0$$ which we prove by induction on $k$. This assertion 
is trivial for $k=0$, so assume it holds for some $k\ge 0$. It follows from $(**)$ and the 
induction hypothesis that
\begin{eqnarray*}
a_{n,k+1} & = & a_{n,k}+a_{n-1,k}\\
          & = & \sum_{l=0}^k{k\choose l}a_{n-l,0}+\sum_{l=0}^k{k\choose l}a_{n-1-l,0}\\
          & = & a_{n,0}+\left(\sum_{l=1}^k{k\choose l}a_{n-l,0}\right)+\left(\sum_{l=1}^k
                {k\choose l-1}a_{n-l,0}\right)+a_{n-1-k,0}\\
          & = & \sum_{l=0}^{k+1}{k+1\choose l}a_{n-l,0}\,,
\end{eqnarray*}
whence the assertion holds also for $k+1$.

(2): In order to prove the assertion, we first derive the identity $$(***)\qquad 
a_{n,l}=\sum_{\nu=0}^n(-1)^\nu a_{n-\nu,l+1}\qquad\forall~n\in\mathbb{N}_0,0\le l
\le k-1$$ from $(**)$ and a telescoping sum argument via the observation that
\begin{eqnarray*}
\sum_{\nu=0}^n(-1)^\nu a_{n-\nu,l+1} & = & \sum_{\nu=0}^n(-1)^\nu a_{n-\nu,l}+
                                           \sum_{\nu=0}^{n-1}(-1)^\nu a_{n-\nu-1,l}\\
                                     & = & a_{n,l}+\sum_{\nu=1}^n(-1)^\nu a_{n-\nu,l}+
                                           \sum_{\mu=1}^n(-1)^{\mu-1} a_{n-\mu,l}\\
                                     & = & a_{n,l}\,.
\end{eqnarray*}
Now we are ready to prove (2). It is a consequence of the slightly more general identity 
$$a_{n,0}=\sum_{\nu=0}^n(-1)^\nu{\nu+l-1\choose l-1}a_{n-\nu,l}\qquad\forall~n\in
\mathbb{N}_0,0\le l\le k$$ which can be shown by finite induction on $l$. For $l=0$ 
the assertion is trivial, so assume it holds for some $0\le l\le k-1$. It follows from 
the induction hypothesis and $(***)$ that
\begin{eqnarray*}
a_{n,0} & = & \sum_{\nu=0}^n(-1)^\nu{\nu+l-1\choose l-1}a_{n-\nu,l}\\
        & = & \sum_{\nu=0}^n(-1)^\nu{\nu+l-1\choose l-1}\left(\sum_{\mu=0}^{n-\nu}
              (-1)^\mu a_{n-\nu-\mu,l+1}\right)\\
        & = & \sum_{\nu=0}^n{\nu+l-1\choose l-1}\left(\sum_{\mu=\nu}^n(-1)^\mu 
              a_{n-\mu,l+1}\right)\\
        & = & \sum_{\mu=0}^n(-1)^\mu\left(\sum_{\nu=0}^\mu{\nu+l-1\choose l-1}\right)
              a_{n-\mu,l+1}\\
        & = & \sum_{\mu=0}^n(-1)^\mu{\mu+l\choose l}a_{n-\mu,l+1}\,, 
\end{eqnarray*}
whence the assertion holds also for $l+1$.\quad $\Box$
\bigskip

We rewrite the first part of Theorem \ref{cohpenv} as an isomorphism theorem similar to a 
factorization theorem of Hochschild and Serre in characteristic zero (see \cite[Theorem 13]{HS}).

\begin{thm}\label{facthm}
Let $L$ be a finite-dimensional Lie algebra over a field $\mathbb{F}$ of prime characteristic 
$p$, let $M$ be a finite-dimensional $L$-module, and let $\mathfrak{L}$ be a $p$-envelope of 
$L$ with $\dim_\mathbb{F}\mathfrak{L}/L<\infty$. Then, for every non-negative integer $n$, 
there is an isomorphism $$H^n(\mathfrak{L},M)\cong\bigoplus_{i+j=n}\Lambda^i(\mathfrak{L}/L)
\otimes_\mathbb{F}H^j(L,M)$$ of $\mathbb{F}$-vector spaces.\quad $\Box$
\end{thm}

\noindent {\it Remark.\/} In particular, Theorem \ref{facthm} shows that $H^n(L,M)$ is
a direct summand of $H^n(\mathfrak{L},M)$ (cf.~\cite[(3.15)]{Leg}).
\bigskip

We conclude this section with some immediate consequences of Theorem \ref{cohpenv} that will 
be useful in Section 5. The first result states that the cohomology of any finite-dimensional
$p$-envelope of a finite-dimensional modular Lie algebra $L$ vanishes if and only if the 
cohomology of $L$ vanishes.

\begin{cor}\label{vancohpenv}
Let $L$ be a finite-dimensional Lie algebra over a field $\mathbb{F}$ of prime characteristic 
$p$, let $M$ be a finite-dimensional $L$-module, and let $\mathfrak{L}$ be a $p$-envelope of 
$L$ with $\dim_\mathbb{F}\mathfrak{L}/L<\infty$. Then $H^j(L,M)=0$ for every $0\le j\le n$ if 
and only if $H^j(\mathfrak{L},M)=0$ for every $0\le j\le n$. In particular, $H^\bullet(L,M)=0$ 
if and only if $H^\bullet(\mathfrak{L},M)=0$.\quad $\Box$
\end{cor}

The low degree cohomology of finite-dimensional $p$-envelopes can be computed as follows (for 
the second part see also \cite[Korollar 3.10]{Leg}).

\begin{cor}\label{lowcohpenv}
Let $L$ be a finite-dimensional Lie algebra over a field $\mathbb{F}$ of prime characteristic 
$p$, let $M$ be a finite-dimensional $L$-module, and let $\mathfrak{L}$ be a $p$-envelope of 
$L$ with $k:=\dim_\mathbb{F}\mathfrak{L}/L<\infty$. Then the following statements hold:
\begin{enumerate}
\item[{\rm (1)}] $H^0(\mathfrak{L},M)\cong H^0(L,M)$.
\item[{\rm (2)}] $\dim_\mathbb{F}H^1(\mathfrak{L},M)=\dim_\mathbb{F}H^1(L,M)+k\dim_\mathbb{F}M^L$.
\item[{\rm (3)}] $\dim_\mathbb{F}H^2(\mathfrak{L},M)=\dim_\mathbb{F}H^2(L,M)+k\dim_\mathbb{F}
                 H^1(L,M)+\frac{k(k-1)}{2}\dim_\mathbb{F}M^L$.\quad $\Box$
\end{enumerate}
\end{cor}


\section{Cohomological Vanishing Theorems}


In this section we reformulate a cohomological vanishing theorem of Dzhumadil'daev by 
using the concept of a universal $p$-envelope. Then we derive some consequences which 
will become useful in the last section. By using the cohomology theory of associative 
algebras Farnsteiner showed that Dzhumadil'daev's vanishing theorem holds more generally 
for arbitrary elements of the center of the augmentation ideal of the universal enveloping 
algebra (see \cite[Theorem 2.1 and Corollary 2.2]{Fa3} as well as \cite{Fa4} for further 
generalizations).

Let $\mathbb{F}$ be a field of prime characteristic $p$. A polynomial of the form $f(t)=
\sum_{n\ge 0}\alpha_nt^{p^n}\in\mathbb{F}[t]$ is called {\it $p$-polynomial\/} (cf.~\cite[
Exercise V.14, p.~196]{Jac}). It is clear from $$\hat{L}=\langle L\rangle_p=\sum_{x\in L}
\sum_{r\in\mathbb{N}_0}\mathbb{F}x^{p^r}$$ (cf.~\cite[Proposition 2.1.3(1)]{SF}) that 
$\hat{L}$ is the $\mathbb{F}$-subspace of $U(L)$ that is generated by the evaluations 
of $p$-polynomials over $\mathbb{F}$ in elements of $L$. This shows that the special
elements in Dzhumadil'daev's cohomological vanishing theorem \cite[Theorem 1]{Dzh1} 
are exactly the elements of the center of the universal $p$-envelope.

\begin{thm}\label{vancoh}
Let $L$ be a Lie algebra over a field of prime characteristic and let $M$ be an $L$-module. 
If there exists an element $z\in C(\hat{L})$ such that $(z)_M$ is invertible, then $H^n(L,M)=0$ 
for every non-negative integer $n$.
\end{thm}

\noindent {\it Proof.\/} Consider a cocycle $\psi\in Z^n(L,M)$ of degree $n$. According to
Lemma \ref{triv}, $\theta^n(z)(\psi)$ is a coboundary. But $z\in C(\hat{L})$ and thus 
$\theta^n(z)(\psi)=(z)_M\circ\psi$. Since $(z)_M$ is invertible and $B^n(L,M)$ is an 
$\hat{L}$-submodule of $C^n(L,M)$, it follows that $\psi\in B^n(L,M)$.\quad $\Box$
\bigskip

In the remainder of this section we apply Theorem \ref{vancoh} to the cohomology of 
simple or indecomposable modules. We begin by proving another result of Dzhumadil'daev 
(cf.~\cite[Theorem 2]{Dzh1}).

\begin{cor}\label{vancohirred}
Let $L$ be a restricted Lie algebra over a field of prime characteristic and let $S$ be 
a simple $L$-module. If $H^n(L,S)\ne 0$ for some non-negative integer, then $S$ is a 
restricted $L$-module.
\end{cor}

\noindent {\it Proof.\/} Suppose that $S$ is not restricted. Then there exists $x\in L$ 
with $(x)_S^p\ne(x^{[p]})_S$. Consider $z:=x^p-x^{[p]}\in C(\hat{L})$. Since $S$ is simple,
it follows from Schur's lemma that $(z)_S$ is invertible. Hence the assertion is a consequence 
of Theorem \ref{vancoh}.\quad $\Box$
\bigskip

Corollary \ref{vancohirred} can be generalized in several directions (see \cite[Theorem 
6.2]{Fa4} or \cite[Corollary 6.4]{Fa4}) but we will need instead an application of the 
following result to $\mathbb{Z}$-graded Lie algebras. 

\begin{cor}\label{vancohadnilp}
Let $L$ be a Lie algebra over a field of prime characteristic and let $M$ be a 
finite-dimensional indecomposable $L$-module. If there exists an ad-nilpotent 
element $x\in L$ such that $(x)_M$ is not nilpotent, then $H^n(L,M)=0$ for every 
non-negative integer $n$.
\end{cor}

\noindent {\it Proof.\/} Since $x$ is ad-nilpotent, there exists a positive integer $q$
such that $(\mathrm{ad}_L\,x)^q=0$, i.e., $z:=x^{p^e}\in C(\hat{L})$ if $p^e\ge q$. 
According to Fitting's lemma (cf.~\cite[p.~37]{Jac}), $M=M_0(z)\oplus M_1(z)$ where 
$$M_0(z):=\bigcup_{n\in\mathbb{N}}\mathrm{Ker}(z)_M^n\qquad\mbox{and}\qquad M_1(z):=
\bigcap_{n\in\mathbb{N}}\mathrm{Im}(z)_M^n\,.$$ It follows from $z\in C(\hat{L})$ that 
$M_0(z)$ and $M_1(z)$ are $L$-submodules of $M$. Since $(x)_M$ is not nilpotent and $M$ 
is a restricted $\hat{L}$-module, $(z)_M$ is also not nilpotent, i.e., $M_0(z)\neq M$. 
But $M$ is indecomposable, so $(z)_M$ is invertible on $M=M_1(z)$; therefore the 
assertion is a consequence of Theorem \ref{vancoh}.\quad $\Box$
\bigskip

\noindent {\it Remark.\/} Corollary \ref{vancohadnilp} holds for arbitrary modules if one 
assumes that $(x)_M$ is invertible (cf.~\cite[Corollary 2 of Theorem 1]{Dzh1}).
\bigskip

\noindent In particular, Corollary \ref{vancohadnilp} applies to finite-dimensional 
$\mathbb{Z}$-graded Lie algebras. Let $$L=\bigoplus_{n\in\mathbb{Z}}L_n$$ be a 
$\mathbb{Z}$-graded Lie algebra and set $$L^+:=\bigoplus_{n>0}L_n\mbox{ and }
L^-:=\bigoplus_{n<0}L_n\,.$$ If $L$ is finite-dimensional, then every element in
$L^+\cup L^-$ is ad-nilpotent so that one can apply Corollary \ref{vancohadnilp} 
in order to prove the following generalization of \cite[Theorem 1.1]{CS}:

\begin{cor}\label{vancohgrad}
Let $L$ be a finite-dimensional $\mathbb{Z}$-graded Lie algebra over a field of prime 
characteristic and let $M$ be a finite-dimensional indecomposable $L$-module. If there 
exists an element $x\in L^+\cup L^-$ such that $(x)_M$ is not nilpotent, then $H^n(L,M)
=0$ for every non-negative integer $n$.\quad $\Box$
\end{cor}

\noindent {\it Remark.\/} Note that a {\it simple\/} $L$-module $S$ is graded if and 
only if $L^+\cup L^-$ acts nilpotently on $S$ (see \cite[Proposition 1.2(1)]{CS}).


\section{Truncated Induced Modules}


In this section we introduce {\it truncated induced modules\/} and {\it truncated 
coinduced modules\/} following \cite{Dzh4}, \cite{Dzh5}, and \cite{FS}. For 
applications of truncated induced modules in the structure theory of modular Lie 
algebras we refer the reader to \cite{Kuz}.

Let $\mathbb{F}$ be a field of prime characteristic $p$ and let $L$ be a finite-dimensional
Lie algebra over $\mathbb{F}$. For every subalgebra $K$ of $L$ with finite cobasis $\{e_1,
\dots,e_k\}$ there exist positive integers $m_1,\dots,m_k$ and $v_i\in U(L)_{(p^{m_i}-1)}$
such that $z_i:=e_i^{p^{m_i}}+v_i\in C(U(L))$ for every $1\le i\le k$. Now consider the 
unital associative subalgebra $\mathcal{O}(L,K)$ of $U(L)$ generated by $K\cup\{z_1,\dots,
z_k\}$.

For $\mathbf{a}\in\mathbb{N}_0^k$ we define $\mathbf{e}^\mathbf{a}:=e_1^{a_1}\cdots 
e_k^{a_k}$ and set $\tau:=(p^{m_1}-1,\dots,p^{m_k}-1)$. It follows from Jacobson's
refinement of the Poincar\'e-Birckhoff-Witt theorem (see \cite[Lemma V.4, p.~189]{Jac}) 
that $U(L)$ is a free left and right $\mathcal{O}(L,K)$-module with basis 
$\{\mathbf{e}^\mathbf{a}\mid 0\le\mathbf{a}\le\tau\}$. As a consequence $$\mathcal{O}
(L,K)\cong\mathbb{F}[z_1,\dots,z_k]\otimes_\mathbb{F}U(K)$$ as unital associative 
$\mathbb{F}$-algebras.

Let $\sigma:K\to\mathbb{F}$ be the Lie algebra homomorphism given by $\sigma(x):=
\mathrm{tr}(\mathrm{ad}_{L/K}\,x)$ for every $x\in K$. Since the mapping $\overline{\sigma}$
defined by $x\mapsto x+\sigma(x)1$ is a Lie algebra homomorphism from $K$ into $U(K)$, 
there exists a unique algebra homomorphism $U(\sigma)$ from $U(K)$ into $U(K)$ that 
extends $\overline{\sigma}$. 

Let $V$ be a $K$-module. Then the action of $U(K)$ on $V$ can be extended to 
$\mathcal{O}(L,K)$ by letting the polynomial algebra $\mathbb{F}[z_1,\dots,z_k]$ 
act via its canonical augmentation mapping. Let $V_\sigma$ be the twisted module 
with $K$-action given by $x\circ v:=x\cdot v+\sigma(x)v$ for every $x\in K$ and 
every $v\in V$. It is clear that $V_\sigma\cong V\otimes_\mathbb{F}\mathbb{F}_\sigma$ 
and that the $\mathcal{O}(L,K)$-action on $V_\sigma$ is given by $\mathrm{id}_{\mathbb{F}
[z_1,\dots,z_k]}\otimes U(\sigma)$.

The following result is a generalization of \cite[Corollary 1.2]{Fe2} to truncated 
induced and coinduced modules.

\begin{thm}\label{frobrec}
Let $L$ be a finite-dimensional Lie algebra over a field of prime characteristic and let 
$K$ be a subalgebra of $L$. Then, for every $K$-module $V$, there are the following 
$L$-module isomorphisms:
\begin{enumerate}
\item[{\rm (1)}] $U(L)\otimes_{\mathcal{O}(L,K)}V\cong\mathrm{Hom}_{\mathcal{O}(L,K)}
                 (U(L),V_{-\sigma})$.
\item[{\rm (2)}] $[U(L)\otimes_{\mathcal{O}(L,K)}V]^*\cong U(L)\otimes_{\mathcal{O}(L,K)}
                 [V^*]_\sigma$.
\item[{\rm (3)}] $[\mathrm{Hom}_{\mathcal{O}(L,K)}(U(L),V)]^*\cong\mathrm{Hom}_{\mathcal{O}
                 (L,K)}(U(L),[V^*]_{-\sigma})$.
\end{enumerate}
\end{thm}

\noindent {\it Proof.\/} (1) is just \cite[Theorem 1.4]{FS} while (2) follows from 
$[U(L)\otimes_{\mathcal{O}(L,K)}V]^*\cong\mathrm{Hom}_{\mathcal{O}(L,K)}(U(L),V^*)$ 
and (1). Finally, (3) is dual to (2).\quad $\Box$
\bigskip

\noindent {\it Remark.\/} Theorem \ref{frobrec}(1) implies that $U(L)\supseteq{\mathcal{O}
(L,K)}$ is a {\it free Frobenius extension of the second kind\/} (cf.~\cite[Theorem 
1.1]{Fe2} for the restricted analogue which is a consequence of Theorem \ref{frobrec}(1) 
and \cite[Proposition 1.5]{FS}).
\bigskip

Let $L$ be a finite-dimensional modular Lie algebra and let $K$ be a subalgebra of 
$L$. Choose a basis $\{e_1,\dots,e_k,e_{k-1},\dots,e_l\}$ of $L$ such that $\{e_{k-1},
\dots,e_l\}$ is a basis of $K$. Then there exist positive integers $m_1,\dots,m_l$ and 
$v_i\in U(L)_{(p^{m_i}-1)}$ such that $z_i:=e_i^{p^{m_i}}+v_i\in C(U(L))$ for every 
$1\le i\le l$. Let $I(L)$ denote the two-sided ideal of $U(L)$ generated by $\{z_1,
\dots,z_l\}$, and set $\overline{U}(L):=U(L)/I(L)$. Similarly, let $I(K)$ denote the 
two-sided ideal of $U(K)$ generated by $\{z_{k+1},\dots,z_l\}$, and set $\overline{U}
(K):=U(K)/I(K)$.
\bigskip

\noindent {\it Remark.\/} Let $L$ be a restricted Lie algebra and let $\chi\in L^*$. 
By choosing $z_i:=e_i^p-e_i^{[p]}-\chi(e_i)^p\cdot 1$ for every $1\le i\le l$, it can 
be seen that $\overline{U}(L)$ and the $\chi$-reduced universal enveloping algebra 
$u(L,\chi)$ are isomorphic as unital associative $\mathbb{F}$-algebras (cf.~\cite[p.~157]{FS}).
\bigskip

Let $V$ be a $K$-module. Then the action of $U(K)$ on $V$ can be extended to $\overline{U}
(K)$ via $z_i\cdot v:=0$ for every $k+1\le i\le l$ and every $v\in V$. Since $I(K)=I(L)
\cap K$, one can consider $\overline{U}(K)$ as a subalgebra of $\overline{U}(L)$ and 
then form induced and coinduced modules called {\it truncated induced modules\/} and 
{\it truncated coinduced modules\/}, respectively. In the special case of restricted 
Lie algebras, the truncated (co)induced modules are just (co)induced modules over reduced 
universal enveloping algebras. The following result (cf.~\cite[p.~35]{Dzh5}) is a 
generalization of \cite[Proposition 1.5]{FS} and can be proved in exactly the same 
manner.

\begin{thm}\label{redtrunc}
Let $L$ be a finite-dimensional Lie algebra over a field of prime characteristic 
and let $K$ be a subalgebra of $L$. Then, for every $K$-module $V$, there are the
following $L$-module isomorphisms:
\begin{enumerate}
\item[{\rm (1)}] $U(L)\otimes_{\mathcal{O}(L,K)}V\cong\overline{U}(L)
                 \otimes_{\overline{U}(K)}V$.
\item[{\rm (2)}] $\mathrm{Hom}_{\mathcal{O}(L,K)}(U(L),V)\cong
                 \mathrm{Hom}_{\overline{U}(K)}(\overline{U}(L),V)$.\quad $\Box$
\end{enumerate}
\end{thm}

\noindent {\it Remark.\/} It is a consequence of Theorem \ref{redtrunc} and Theorem 
\ref{frobrec}(1) that $\overline{U}(L)\supseteq\overline{U}(K)$ is a {\it free 
Frobenius extension of the second kind\/}.
\bigskip

In the last section we will need Shapiro's lemma for truncated induced modules (see 
\cite[Theorem 2.1]{FS} or \cite[Theorem in \S5 and Corollary 1 in \S3]{Dzh4}):

\begin{thm}\label{shapiro}
Let $L$ be a finite-dimensional Lie algebra over a field of prime characteristic, let 
$K$ be a subalgebra of $L$, and let $V$ be a $K$-module. Then, for every non-negative 
integer $n$, there is an isomorphism $$H^n(L,U(L)\otimes_{\mathcal{O}(L,K)}V)\cong 
\bigoplus_{i+j=n}\Lambda^i(L/K)\otimes_\mathbb{F}H^j(K,V_{-\sigma})$$ of $\mathbb{F}$-vector 
spaces.
\end{thm}

\noindent {\it Proof.\/} Let $n$ be a non-negative integer. Then the following isomorphisms
can be obtained from Theorem \ref{frobrec}(1) in conjunction with Shapiro's lemma and the 
K\"unneth formula (cf.~\cite[Theorem 3.6.3]{Wei}): 
\begin{eqnarray*}
H^n(L,U(L)\otimes_{\mathcal{O}(L,K)}V) & \cong & \mathrm{Ext}_{U(L)}^n(\mathbb{F},U(L)
                                                 \otimes_{\mathcal{O}(L,K)}V)\\
                                       & \cong & \mathrm{Ext}_{U(L)}^n(\mathbb{F},
                                                 \mathrm{Hom}_{\mathcal{O}(L,K)}
                                                 (U(L),V_{-\sigma}))\\
                                       & \cong & \mathrm{Ext}_{\mathcal{O}(L,K)}^n
                                                 (\mathbb{F},V_{-\sigma})\\
                                       & \cong & \mathrm{Ext}_{\mathbb{F}[z_1,\dots,z_k]
                                                 \otimes_\mathbb{F}U(K)}^n(\mathbb{F}
                                                 \otimes_\mathbb{F}V_{-\sigma})\\
                                       & \cong & \bigoplus_{i+j=n}\mathrm{Ext}_{\mathbb{F}
                                                 [z_1,\dots,z_k]}^i(\mathbb{F},\mathbb{F})
                                                 \otimes_\mathbb{F}\mathrm{Ext}_{U(K)}^j
                                                 (\mathbb{F},V_{-\sigma})\\
                                       & \cong & \bigoplus_{i+j=n}\Lambda^i(L/K)
                                                 \otimes_\mathbb{F}H^j(K,V_{-\sigma})\,.
                                                 \quad\Box
\end{eqnarray*}
\smallskip

We conclude this section by considering certain truncated induced modules over the Zassenhaus 
algebras. Let $\mathcal{O}(W(m),\mathfrak{b}(m))$ denote the unital associative subalgebra of 
$U(W(m))$ generated by $\mathfrak{b}(m)$ and $z_{-1}=e_{-1}^{p^m}\in C(\hat{W}(m))$. Moreover, 
for every $\lambda\in\mathbb{F}$, let $F_\lambda:=\mathbb{F}1_\lambda$ be the one-dimensional 
$\mathfrak{b}(m)$-module defined by $e_0\cdot 1_\lambda:=\lambda 1_\lambda$ and by $e_i\cdot 
1_\lambda:=0$ for every $1\le i\le p^m-2$. Then $F_\lambda$ is a unital $\mathcal{O}(W(m),
\mathfrak{b}(m))$-module via $z_{-1}\cdot 1_\lambda:=0$. In analogy with the baby Verma 
modules of Lie algebras of reductive groups, the truncated induced module $V(\lambda):=U(W(m))
\otimes_{\mathcal{O}(W(m),\mathfrak{b}(m))}F_\lambda$ is called a {\it generalized baby Verma 
module\/}.

The next result will be useful in Section 5 and is included here for the convenience
of the reader. Recall that the divided power algebra $A(m)$ is a natural $W(m)$-module
via $(f\partial)\cdot g:=f\partial(g)$ for all $f,g\in A(m)$. 

\begin{pro}\label{adjnat}
Let $\mathbb{F}$ be a field of prime characteristic $p$ and let $m$ be a positive integer. 
Then there are the following $W(m)$-module isomorphisms:
\begin{enumerate}
\item[{\rm (1)}] $W(m)\cong V(p-2)$.
\item[{\rm (2)}] $W(m)^*\cong V(1)$.
\item[{\rm (3)}] $A(m)\cong V(p-1)$.
\item[{\rm (4)}] $A(m)^*\cong V(0)$.
\end{enumerate}
\end{pro}

\noindent {\it Proof.\/} Set $W:=W(m)$, $\mathfrak{b}:=\mathfrak{b}(m)$, and $A:=A(m)$.

(1): Consider the $\mathcal{O}(W,\mathfrak{b})$-module homomorphism $\varphi:F_{p-2}
\to W$ given by $1_{p-2}\mapsto e_{p^m-2}$. Then the universal property of induced 
modules (cf.~\cite[(I.3.3)]{Fe1}) yields a $W$-module homomorphism $\Phi:V(p-2)\to 
W$ such that $\Phi(1\otimes 1_{p-2})=\varphi(1_{p-2})=e_{p^m-2}$. In particular, 
$\Phi\neq 0$. Since $W$ is simple, $\Phi$ is surjective. But $\dim_\mathbb{F}V(p-2)
=p^m=\dim_\mathbb{F}W$ and therefore $\Phi$ is bijective. (By employing Lemma 
\ref{ind} below, the assertion is also an immediate consequence of $\Phi\neq 0$ 
and Schur's lemma.)

(3): Consider the $\mathcal{O}(W,\mathfrak{b})$-module homomorphism $\psi:F_{p-1}
\to A$ given by $1_{p-1}\mapsto x^{(p^m-1)}$. Then the universal property of induced 
modules (cf.~\cite[(I.3.3)]{Fe1}) yields a $W$-module homomorphism $\Psi:V(p-1)\to 
A$ such that $\Psi(1\otimes 1_{p-1})=\psi(1_{p-1})=x^{(p^m-1)}$. In particular, 
$\Psi\neq 0$. Since $$\Psi(e_{-1}^a\otimes 1_{p-1})=e_{-1}^a\cdot x^{(p^m-1)}=
x^{(p^m-1-a)}$$ for every $0\le a\le p^m-1$, $\Psi$ is surjective. But 
$\dim_\mathbb{F}V(p-1)=p^m=\dim_\mathbb{F}A$ and thus $\Psi$ is bijective.

(2) \& (4): Finally, consider the Lie algebra homomorphism $\sigma:\mathfrak{b}\to
\mathbb{F}$ defined by $\sigma(x):=\mathrm{tr}(\mathrm{ad}_{W/\mathfrak{b}}\,x)$ for 
every $x\in\mathfrak{b}$. We have $[e_i,e_{-1}]=-e_{i-1}$ and therefore $\sigma(e_i)
=-\delta_{i0}$ for every $0\le i\le p^m-2$. Hence $[F_{-\lambda}]_\sigma=F_{p-1-\lambda}$ 
and thus Theorem \ref{frobrec}(2) yields $V(\lambda)^*\cong V(p-1-\lambda)$ for every 
$\lambda\in\mathbb{F}_p$. Then (2) and (4) are immediate consequences of this isomorphism 
in conjunction with (1) and (3), respectively.\quad $\Box$
\bigskip

\noindent {\it Remark.\/} It follows from the first two parts of Proposition \ref{adjnat}
that $W(m)^*\cong W(m)$ if and only if $p=3$, i.e., Zassenhaus algebras have
non-degenerate invariant bilinear forms if and only if $p=3$ (cf.~\cite[Theorem 1]{Dzh2}
or \cite[Theorem 4.6.3]{SF}).
\bigskip

It is well-known that the generalized baby Verma modules of a Zassenhaus algebra $W(m)$ 
can be realized on the divided power algebra $A(m)$ (see \cite[p.~134]{Dzh1}). Let 
$\lambda\in\mathbb{F}$ and let $A_\lambda(m)$ denote the vector space $A(m)$ with 
$W(m)$-action defined by $(f\partial)\cdot g:=f\partial(g)+\lambda\partial(f)g$ for 
all $f,g\in A(m)$. As in the proof of Proposition \ref{adjnat}(3), one obtains the
following isomorphisms:

\begin{pro}\label{divpowalg}
Let $\mathbb{F}$ be a field of prime characteristic $p$ and let $m$ be a positive integer. 
Then, for every $\lambda\in\mathbb{F}$, there is a $W(m)$-module isomorphism $V(\lambda)
\cong A_{\lambda+1}(m)$.\quad $\Box$
\end{pro}

Finally, it is easy to see that the generalized baby Verma modules of Zassenhaus algebras 
are induced modules for restricted universal enveloping algebras of their minimal 
$p$-envelopes. Note that $F_\lambda$ is a restricted $\mathfrak{b}(m)$-module if and 
only if $\lambda\in\mathbb{F}_p$.

\begin{pro}\label{verma}
Let $\mathbb{F}$ be a field of prime characteristic $p$ and let $m$ be a positive integer.
Then, for every $\lambda\in\mathbb{F}_p$, there is a $W(m)$-module isomorphism $V(\lambda)
\cong u(\mathfrak{W}(m))\otimes_{u(\mathfrak{b}(m))}F_\lambda$ where $z_{-1}:=e_{-1}^{p^m}$.
\quad $\Box$
\end{pro}

\noindent {\it Remark.\/} In general, there is a similar relationship between truncated 
(co)induced modules of modular Lie algebras $L\supseteq K$ and (co)induced modules over 
restricted universal enveloping algebras of certain $p$-envelopes $\mathfrak{L}\supseteq
\mathfrak{K}$ of $L$ and $K$, respectively. Namely, Farnsteiner \cite[Theorem 2.3]{Fa5} 
proved that for every finite-dimensional $K$-module $V$ there exist finite-dimensional 
$p$-envelopes $\mathfrak{K}\subseteq\mathfrak{L}$ of $K$ and $L$, respectively, such that 
$V$ is a restricted $K$-module and $U(L)\otimes_{\mathcal{O}(L,K)}V\cong u(\mathfrak{L})
\otimes_{u(\mathfrak{K})}V$ as $L$-modules.


\section{Cohomology of Zassenhaus Algebras}


In this section we compute the cohomology of Zassenhaus algebras with coefficients 
in generalized baby Verma modules and in simple modules.

The first result reduces the computation of the cohomology of the Zassenhaus algebra $W(m)$
with coefficients in a generalized baby Verma module to the $\mathfrak{t}$-invariants of 
the cohomology of the maximal nilpotent subalgebra $\mathfrak{u}(m)$ with coefficients 
in the one-dimensional trivial module (cf.~\cite[Theorem 5]{Dzh1}).

\begin{thm}\label{cohzas}
Let $\mathbb{F}$ be a field of prime characteristic $p$, let $m$ be a positive integer, 
and let $\lambda\in\mathbb{F}$. Then, for every non-negative integer $n$, there is an 
isomorphism 
\begin{eqnarray*}
H^n(W(m),V(\lambda)) & \cong & H^n(\mathfrak{b}(m),F_{\lambda+1})\oplus H^{n-1}
                               (\mathfrak{b}(m),F_{\lambda+1})\\
                     & \cong & [H^n(\mathfrak{u}(m),F_{\lambda+1})\oplus H^{n-1}
                               (\mathfrak{u}(m),F_{\lambda+1})^{\oplus 2}\oplus 
                               H^{n-2}(\mathfrak{u}(m),F_{\lambda+1})]^\mathfrak{t}
\end{eqnarray*} 
of $\mathbb{F}$-vector spaces.
\end{thm}

\noindent {\it Proof.\/} Set $W:=W(m)$, $\mathfrak{b}:=\mathfrak{b}(m)$, and $\mathfrak{u}
:=\mathfrak{u}(m)$. Let $\sigma:\mathfrak{b}\to\mathbb{F}$ be defined by $\sigma(x):=
\mathrm{tr}(\mathrm{ad}_{W/\mathfrak{b}}\,x)$ for every $x\in\mathfrak{b}$. As in the 
proof of Proposition \ref{adjnat}, we have $\sigma(e_i)=-\delta_{i0}$ for every $0\le i
\le p^m-2$. Hence $[F_\lambda]_{-\sigma}=F_{\lambda+1}$ and the first isomorphism follows 
from Theorem \ref{shapiro}. Moreover, the second isomorphism can be obtained from Lemma 
\ref{codim1} and the fact that $\mathfrak{u}H^n(\mathfrak{u},F_{\lambda+1})=0$.\quad 
$\Box$
\bigskip

Later on in this section we will need the invariant spaces of generalized baby Verma 
modules for the Zassenhaus algebras and their minimal $p$-envelopes:

\begin{pro}\label{invzas}
Let $\mathbb{F}$ be a field of prime characteristic $p$, let $m$ be a positive integer, 
and let $\lambda\in\mathbb{F}$. Then there are isomorphisms
\begin{eqnarray*}
V(\lambda)^{\mathfrak{W}(m)}\cong V(\lambda)^{W(m)}\cong
\left\{
\begin{array}{cl}
\mathbb{F} & \mbox{if }\lambda=p-1\\ 
0 & \mbox{if }\lambda\neq p-1\,.
\end{array}
\right.
\end{eqnarray*}
of $\mathbb{F}$-vector spaces.
\end{pro}

\noindent {\it Proof.\/} It is a consequence of Corollary \ref{lowcohpenv}(1) and Theorem 
\ref{cohzas} that
\begin{eqnarray*}
V(\lambda)^{\mathfrak{W}(m)}\cong V(\lambda)^{W(m)}\cong H^0(W(m),V(\lambda))\cong 
H^0(\mathfrak{u}(m),F_{\lambda+1})^\mathfrak{t}
\cong(F_{\lambda+1})^\mathfrak{t}\cong
\left\{
\begin{array}{cl}
\mathbb{F} & \mbox{if }\lambda=p-1\\ 
0 & \mbox{if }\lambda\neq p-1\,.\quad\Box
\end{array}
\right.
\end{eqnarray*}
\smallskip

For $p>3$ the following result was proved by Dzhumadil'daev (see \cite[Corollary 
2 of Theorem 5]{Dzh1}) while for $p=3$ it is contained in \cite[Satz 5.17]{Leg}.

\begin{thm}\label{1cohzas}
Let $\mathbb{F}$ be a field of prime characteristic $p>2$, let $m$ be a positive integer, 
and let $\lambda\in\mathbb{F}$. Then the following statements hold:
\begin{enumerate}
\item[{\rm (1)}] If $p=3$, then 
                 \begin{eqnarray*}
                 \dim_\mathbb{F}H^1(W(m),V(\lambda))=
                 \left\{
                 \begin{array}{cl}
                 1 & \mbox{if }\lambda=0\\ 
                 m-1 & \mbox{if }\lambda=1\\
                 2 & \mbox{if }\lambda=2\\
                 0 & \mbox{if }\lambda\ne 0,1,2\,.
                 \end{array}
                 \right.
                 \end{eqnarray*}
\item[{\rm (2)}] If $p>3$, then
                 \begin{eqnarray*}
                 \dim_\mathbb{F}H^1(W(m),V(\lambda))=
                 \left\{
                 \begin{array}{cl}
                 1 & \mbox{if }\lambda=0\mbox{ or }\lambda=1\\ 
                 m-1 & \mbox{if }\lambda=p-2\\
                 2 & \mbox{if }\lambda=p-1\\
                 0 & \mbox{if }\lambda\ne 0,1,p-2,p-1\,.
                 \end{array}
                 \right.
                 \end{eqnarray*}
\end{enumerate}
\end{thm}

\noindent {\it Proof.\/} Set $W:=W(m)$ and $\mathfrak{u}:=\mathfrak{u}(m)$. Since $\mathfrak{u}$ 
acts trivially on $F_{\lambda+1}$, it follows from Theorem \ref{cohzas} that
\begin{eqnarray*}
H^1(W,V(\lambda)) & \cong & H^1(\mathfrak{u},F_{\lambda+1})^\mathfrak{t}\oplus 
                            H^0(\mathfrak{u},F_{\lambda+1})^\mathfrak{t}\oplus 
                            H^0(\mathfrak{u},F_{\lambda+1})^\mathfrak{t}\\ 
                  & \cong & H^1(\mathfrak{u},F_{\lambda+1})^\mathfrak{t}\oplus 
                            (F_{\lambda+1})^\mathfrak{t}\oplus (F_{\lambda+1})^\mathfrak{t}\,.
\end{eqnarray*} 
For the last two summands we have that
\begin{eqnarray*}
(+)\qquad (F_{\lambda+1})^\mathfrak{t}\cong
\left\{
\begin{array}{cl}
\mathbb{F} & \mbox{if }\lambda=p-1\\ 
0 & \mbox{if }\lambda\ne p-1\,. 
\end{array}
\right.
\end{eqnarray*}
Since $\mathfrak{u}$ acts trivially on $F_{\lambda+1}$, we obtain for the first summand that 
$$H^1(\mathfrak{u},F_{\lambda+1})^\mathfrak{t}\cong\mathrm{Hom}_\mathbb{F}(\mathfrak{u}/
[\mathfrak{u},\mathfrak{u}],F_{\lambda+1})^\mathfrak{t}\,,$$ where $\mathfrak{t}$ acts via
$(e_0\cdot\eta)(\overline{x})=e_0\cdot\eta(\overline{x})-\eta(\overline{[e_0,x]})$ for every 
$x\in\mathfrak{u}$. (Here $\overline{x}$ denotes the coset of $x$ in $\mathfrak{u}/[\mathfrak{u},
\mathfrak{u}]$.) By a straightforward computation (see \cite[Lemma 5.14]{Leg}) one can show that 
\begin{eqnarray*}
\mathfrak{u}/[\mathfrak{u},\mathfrak{u}]\cong
\left\{
\begin{array}{cl}
\mathbb{F}\overline{e_1}\oplus\bigoplus_{r=1}^{m-1}\mathbb{F}\overline{e_{3^r-1}} & 
\mbox{if }p=3\\
\mathbb{F}\overline{e_1}\oplus\mathbb{F}\overline{e_2}\oplus\bigoplus_{r=1}^{m-1}\mathbb{F}
\overline{e_{p^r-1}} & \mbox{if }p>3\,.
\end{array}
\right.
\end{eqnarray*}
Suppose now that $\eta\in\mathrm{Hom}_\mathbb{F}(\mathfrak{u}/[\mathfrak{u},\mathfrak{u}],
F_{\lambda+1})^\mathfrak{t}$, i.e., $e_0\cdot\eta=0$. Then $$0=(e_0\cdot\eta)(\overline{e_i})
=e_0\cdot\eta(\overline{e_i})-\eta(\overline{[e_0,e_i]})=(\lambda+1-i)\eta(\overline{e_i})$$ 
for $i\in\{1,2,p^r-1\mid 1\le r\le m-1\}$. In particular, $\mathrm{Hom}_\mathbb{F}(\mathfrak{u}/
[\mathfrak{u},\mathfrak{u}],F_{\lambda+1})^\mathfrak{t}=0$ unless $\lambda\in\mathbb{F}_p$.
Moreover, we obtain the following results which, in conjunction with $(+)$, finish the 
proof of the theorem.
\smallskip

\noindent (1) If $p=3$, then 
\begin{eqnarray*}
\dim_\mathbb{F}\mathrm{Hom}_\mathbb{F}(\mathfrak{u}/[\mathfrak{u},\mathfrak{u}],
F_{\lambda+1})^\mathfrak{t}=
\left\{
\begin{array}{cl}
1 & \mbox{if }\lambda=0\\ 
m-1 & \mbox{if }\lambda=1\,.\\
0 & \mbox{if }\lambda\ne 0,1
\end{array}
\right.
\end{eqnarray*}
(2) If $p>3$, then
\begin{eqnarray*}
\dim_\mathbb{F}\mathrm{Hom}_\mathbb{F}(\mathfrak{u}/[\mathfrak{u},\mathfrak{u}],
F_{\lambda+1})^\mathfrak{t}=
\left\{
\begin{array}{cl}
1 & \mbox{if }\lambda=0\mbox{ or }\lambda = 1\\ 
m-1 & \mbox{if }\lambda=p-2\,.\hspace{1cm}\Box\\
0 & \mbox{if }\lambda\ne 0,1,p-2
\end{array}
\right.
\end{eqnarray*}
\medskip

\noindent {\it Remark.\/} It is an immediate consequence of Theorem \ref{1cohzas} 
and Proposition \ref{adjnat}(1) that $$\dim_\mathbb{F}H^1(W(m),W(m))=m-1\,.$$ In 
fact, the outer derivations of $W(m)$ are induced by $(\mathrm{ad}\,e_{-1})^{p^r}$ 
for any $1\le r\le m-1$ (see also the remarks at the end of Section 1).
\bigskip

In view of the discussion before Lemma \ref{triv}, every $W(m)$-module is a 
$\mathfrak{W}(m)$-module for the minimal $p$-envelope $\mathfrak{W}(m)$ of $W(m)$. 
In particular, $V(\lambda)$ is a $\mathfrak{W}(m)$-module for every $\lambda\in
\mathbb{F}$. The next result is a consequence of Corollary \ref{lowcohpenv}(2),
Theorem \ref{1cohzas}, and Proposition \ref{invzas}.

\begin{thm}\label{1cohreszas}
Let $\mathbb{F}$ be a field of prime characteristic $p>2$, let $m$ be a positive integer, 
and let $\lambda\in\mathbb{F}$. Then the following statements hold:
\begin{enumerate}
\item[{\rm (1)}] If $p=3$, then 
                 \begin{eqnarray*}
                 \dim_\mathbb{F}H^1(\mathfrak{W}(m),V(\lambda))=
                 \left\{
                 \begin{array}{cl}
                 1 & \mbox{if }\lambda=0\\ 
                 m-1 & \mbox{if }\lambda=1\\
                 m+1 & \mbox{if }\lambda=2\\
                 0 & \mbox{if }\lambda\ne 0,1,2\,.
                 \end{array}
                 \right.
                 \end{eqnarray*}
\item[{\rm (2)}] If $p>3$, then
                 \begin{eqnarray*}
                 \dim_\mathbb{F}H^1(\mathfrak{W}(m),V(\lambda))=
                 \left\{
                 \begin{array}{cl}
                 1 & \mbox{if }\lambda=0\mbox{ or }\lambda=1\\ 
                 m-1 & \mbox{if }\lambda=p-2\\
                 m+1 & \mbox{if }\lambda=p-1\\
                 0 & \mbox{if }\lambda\ne 0,1,p-2,p-1\,.\quad\Box
                 \end{array}
                 \right.
                 \end{eqnarray*}
\end{enumerate}
\end{thm}

\noindent {\it Remark.\/} According to Proposition \ref{verma}, $V(\lambda)\cong 
u(\mathfrak{W}(m))\otimes_{u(\mathfrak{b}(m))}F_\lambda$ is a restricted baby 
Verma module for $\mathfrak{W}(m)$ if $\lambda\in\mathbb{F}_p$. This isomorphism, 
in conjunction with \cite[Proposition 1.5]{FS} and Theorem \ref{shapiro}, can be 
employed to give a proof of Theorem \ref{1cohreszas} for $\lambda\in\mathbb{F}_p$ 
which does neither use Theorem \ref{1cohzas} nor the results in Section 2.
\bigskip 

Following Hochschild \cite{Ho} we define the {\it restricted cohomology\/} of a 
restricted Lie algebra $L$ with coefficients in a restricted $L$-module $M$ by means 
of $$H_*^n(L,M):=\mathrm{Ext}_{u(L)}^n(\mathbb{F},M)\qquad\forall~n\in\mathbb{N}_0\,.$$
Observe that $V(\lambda)$ is a restricted $\mathfrak{W}(m)$-module if and only if 
$\lambda\in\mathbb{F}_p$. Then Theorem \ref{1cohreszas} can be used to determine 
the dimensions of the restricted $1$-cohomology of the minimal $p$-envelopes of the 
Zassenhaus algebras with coefficients in restricted baby Verma modules (cf.~\cite[Satz 
III.3.3]{Fe1} for $m=1$ and $p>3$).

\begin{cor}\label{1rescohreszas}
Let $\mathbb{F}$ be a field of prime characteristic $p>2$, let $m$ be a positive integer, 
and let $\lambda\in\mathbb{F}_p$. Then the following statements hold:
\begin{enumerate}
\item[{\rm (1)}] If $p=3$, then
                 \begin{eqnarray*}
                 \dim_\mathbb{F}H_*^1(\mathfrak{W}(m),V(\lambda))=
                 \left\{
                 \begin{array}{cl}
                 1 & \mbox{if }\lambda=0\\
                 m-1 & \mbox{if }\lambda=1\\
                 0 & \mbox{if }\lambda=2\,.
                 \end{array}
                 \right.
                 \end{eqnarray*}
\item[{\rm (2)}] If $p>3$, then
                 \begin{eqnarray*}
                 \dim_\mathbb{F}H_*^1(\mathfrak{W}(m),V(\lambda))=
                 \left\{
                 \begin{array}{cl}
                 1 & \mbox{if }\lambda=0\mbox{ or }\lambda=1\\ 
                 m-1 & \mbox{if }\lambda=p-2\\
                 0 & \mbox{if }\lambda\neq 0,1,p-2\,.
                 \end{array}
                 \right.
                 \end{eqnarray*}
\end{enumerate}
\end{cor} 

\noindent {\it Proof.\/} Set $\mathfrak{W}:=\mathfrak{W}(m)$, $\mathfrak{b}:=\mathfrak{b}
(m)$, and $\mathfrak{u}:=\mathfrak{u}(m)$. According to Proposition \ref{verma}, $V(\lambda)
\cong u(\mathfrak{W})\otimes_{u(\mathfrak{b})}F_\lambda$ is a restricted baby Verma module 
for $\mathfrak{W}$ if $\lambda\in\mathbb{F}_p$. Then, from Shapiro's lemma (cf.~\cite[Corollary 
1.4]{Fe2}) and the Hochschild-Serre spectral sequence for restricted cohomology in conjunction 
with \cite[Corollary 3.6]{Fe2}, \cite[Proposition 2.7]{Fe2}, and the definition of the 
$p$-mapping of $\mathfrak{u}$, one obtains that 
\begin{eqnarray*}
H_*^1(\mathfrak{W},V(\lambda)) & \cong & H_*^1(\mathfrak{W},u(\mathfrak{W})\otimes_{u(\mathfrak{b})}
                                         F_\lambda)\\
                               & \cong & H_*^1(\mathfrak{b},F_{\lambda+1})\\
                               & \cong & H_*^1(\mathfrak{u},F_{\lambda+1})^\mathfrak{t}\\
                               & \cong & \mathrm{Hom}_\mathfrak{t}(\mathfrak{u}/[\mathfrak{u},\mathfrak{u}]
                                         +\langle\mathfrak{u}^{[p]}\rangle_\mathbb{F},F_{\lambda+1})\\
                               & \cong & \mathrm{Hom}_\mathfrak{t}(\mathfrak{u}/[\mathfrak{u},\mathfrak{u}],
                                         F_{\lambda+1})\,. 
\end{eqnarray*}
Now the assertions follow from the dimensions formulas at the end of the proof of 
Theorem \ref{1cohzas}.\quad $\Box$
\smallskip

\noindent {\it Remark.\/} For $\lambda\neq p-1$, Corollary \ref{1rescohreszas} is an 
immediate consequence of Theorem \ref{1cohreszas} and Proposition \ref{invzas} in 
conjunction with the six-term exact sequence of Hochschild \cite[p.~575]{Ho} or the 
Appendix of \cite{Fe3}.
\bigskip

Let $S$ be a simple $W(m)$-module. According to Corollary \ref{vancohgrad}, $H^\bullet
(W(m),S)=0$ unless $e_{-1}$ and $\mathfrak{u}(m)$ act nilpotently on $S$. In the 
following we determine the simple modules of the Zassenhaus algebras on which $e_{-1}$ 
and $\mathfrak{u}(m)$ act nilpotently (cf.~\cite[Lemma 5.2]{Leg}). 

\begin{lem}\label{highwt}
Let $\mathbb{F}$ be an algebraically closed field of prime characteristic $p$ and let 
$S$ be a simple $W(m)$-module on which $e_{-1}$ and $\mathfrak{u}(m)$ act nilpotently. 
Then there is a non-zero element $v\in S$ and there is an element $\lambda\in\mathbb{F}$ 
such that $e_0\cdot v=\lambda v$ and $e_i\cdot v=0$ for every $1\le i\le p^m-2$.
\end{lem}

\noindent {\it Proof.\/} Set $W:=W(m)$ and $\mathfrak{u}:=\mathfrak{u}(m)$. Because $S$ 
is simple, there exists a non-zero element $w\in S$ such that $S=U(W)w$. According to 
\cite[Theorem 3.4.3]{SF}, the $\mathbb{Z}$-grading of $W$ induces a $\mathbb{Z}$-grading 
on $U(W)$. Since $\mathfrak{u}$ acts nilpotently on $S$, the integer $a_*:=\min\{a\in
\mathbb{N}\mid e_i^a\cdot S=0\mbox{ for every }1\le i\le p^m-2\}$ exists. Then $U(W)_n
\cdot w=0$ if $n\ge a_*(p^m-1)^2$. Namely, if $e_{-1}^{a_{-1}}e_0^{a_0}\cdots 
e_{p^m-2}^{a_{p^m-2}}$ is a basis element of $U(W)_n$ with $a_i<a_*$ for every $1\le 
i\le p^m-2$, then $$n=\sum_{i=-1}^{p^m-2}ia_i\le\sum_{i=1}^{p^m-2}ia_i<a_*\sum_{i=1}^{p^m-2}
i=\frac{(p^m-2)(p^m-1)}{2}a_*<a_*(p^m-1)^2\,.$$ Now set $n_*:=\max\{n\in\mathbb{Z}\mid 
U(W)_n\cdot w\ne 0\}$, which exists by the above, and consider $M:=U(W)_{n_*}w$. Then 
$e_0\cdot M\subseteq M$. Since $S$ is finite-dimensional (cf.~\cite[Theorem 5.2.4]{SF}), 
$M$ is also finite-dimensional. But $\mathbb{F}$ is algebraically closed, so $e_0$ has 
an eigenvector $v\in M$ with eigenvalue $\lambda\in\mathbb{F}$. Finally, $e_i\cdot v\in 
U(W)_{n_*+i}\cdot w=0$ for every $1\le i\le p^m-2$.\quad $\Box$
\bigskip

\noindent {\it Remark.\/} Note that a simple $W(m)$-module $S$ over an algebraically closed 
ground field is graded if and only if $e_{-1}$ and $\mathfrak{u}(m)$ act nilpotently on $S$ 
(cf.~\cite[Proposition 1.2(1)]{CS}). 
\bigskip

Let $S$ be a simple $W(m)$-module as in Lemma \ref{highwt}. Consider $z_0:=e_0^p-e_0\in
\hat{W}(m)$. Then $[z_0,e_i]=(i^p-i)e_i=0$ for every $-1\le i\le p^m-2$ and therefore 
$z_0\in C(\hat{W}(m))$. Moreover, $z_0\cdot v=(\lambda^p-\lambda)v$. By virtue of Theorem 
\ref{vancoh}, $H^\bullet(W(m),S)=0$ unless $\lambda^p=\lambda$, i.e., $\lambda$ belongs
to the prime field $\mathbb{F}_p$ (cf.~\cite[Lemma 5.3]{Leg}). Finally, consider $z_{-1}
:=e_{-1}^{p^m}\in C(\hat{W}(m))$ (cf.~Section 1). Since $S$ is simple, it follows from 
Schur's lemma and the remark after Corollary \ref{vancohadnilp} that $H^\bullet(W(m),S)
=0$ unless $e_{-1}^{p^m}\cdot v=0$ (cf.~\cite[Lemma 5.4]{Leg}).

In the following let $S(\lambda)$ denote the simple $W(m)$-module with a non-zero element 
$v\in S(\lambda)$ such that $$e_{-1}^{p^m}\cdot v=0\,,$$ $$e_0\cdot v=\lambda v\mbox{ for 
some }\lambda\in\mathbb{F}_p\,,\mbox{ and}$$ $$e_i\cdot v=0\mbox{ for every }1\le i\le 
p^m-2\,.$$ Let $r:=\max\{n\in\mathbb{N}_0\mid e_{-1}^n\cdot v\ne 0\}\le p^m-1$ and consider 
$$M:=\sum_{n=0}^r\mathbb{F}e_{-1}^n\cdot v\,.$$ It then follows from the Cartan-Weyl identity 
(cf.~\cite[Proposition 1.1.3(4)]{SF}) that 
\begin{eqnarray*}
e_i\cdot(e_{-1}^n\cdot v) & = & \sum_{j=0}^n(-1)^{n-j}{n\choose j}e_{-1}^j
                                (\mathrm{ad}\,e_{-1})^{n-j}(e_i)\cdot v\\
                          & = & \sum_{j=0}^n(-1)^{n-j}{n\choose j}(e_{-1}^je_{i-n+j})\cdot v\\
                          & = & (-1)^{i+1}{n\choose n-1-i}(e_{-1}^{n-1-i}e_{-1})(v)\\
                          &   & +(-1)^i{n\choose n-i}(e_{-1}^{n-i}e_0)\cdot v\\
                          & = & (-1)^i[\lambda{n\choose i}-{n\choose i+1}]e_{-1}^{n-i}\cdot v 
\end{eqnarray*}
for every $0\le i\le p^m-2$ (cf.~\cite[Lemma III.3.1]{Fe1} for $m=1$ and \cite[Lemma 5.5]{Leg}
for the general case). Consequently, $M$ is a non-zero submodule of $S(\lambda)$ and as the 
latter is simple, $M=S(\lambda)$.

Most of the generalized baby Verma modules of Zassenhaus algebras are simple (cf.~\cite[Lemma 
III.3.2a)]{Fe1} for $m=1$ and \cite[Lemma 5.6]{Leg} for the general case).

\begin{lem}\label{ind}
If $\lambda\not\in\{0,p-1\}$, then $V(\lambda)$ is a simple $W(m)$-module.
\end{lem}

\noindent {\it Proof.\/} Set $W:=W(m)$ and $\mathfrak{b}:=\mathfrak{b}(m)$. Then 
$U(W)$ is a free right $\mathcal{O}(W,\mathfrak{b})$-module with basis $\{e_{-1}^n
\mid 0\le n\le p^m-1\}$. Hence $\{e_{-1}^n\otimes 1_\lambda\mid 0\le n\le p^m-1\}$ 
is a basis of $V(\lambda)$ over $\mathbb{F}$. As above it follows from the Cartan-Weyl 
identity that $$e_ie_{-1}^n\otimes 1_\lambda=(-1)^i[\lambda{n\choose i}-{n\choose 
i+1}]e_{-1}^{n-i}\otimes 1_\lambda$$ for every $0\le i\le p^m-2$.

Let $M$ be a non-zero $W$-submodule of $V(\lambda)$ and let $\sum_{n=s}^{p^m-1}
\alpha_ne_{-1}^n\otimes 1_\lambda\in M$ with $\alpha_s\ne 0$. Then $$\alpha_s
e_{-1}^{p^m-1}\otimes 1_\lambda=\sum_{n=s}^{p^m-1}\alpha_ne_{-1}^{n+p^m-1-s}
\otimes 1_\lambda=e_{-1}^{p^m-1-s}\cdot\left(\sum_{n=s}^{p^m-1}\alpha_ne_{-1}^n
\otimes 1_\lambda\right)\in M\,.$$ But $\alpha_s\ne 0$, so $e_{-1}^{p^m-1}\otimes 
1_\lambda\in M$.

Now, for every $0\le i\le p^m-2$, consider
\begin{eqnarray*}
e_ie_{-1}^{p^m-1}\otimes 1_\lambda & = & (-1)^i[\lambda{p^m-1\choose i}-{p^m-1\choose i+1}]
                                         e_{-1}^{p^m-1-i}\otimes 1_\lambda\\ 
                                   & = & (-1)^i[\lambda(-1)^i-(-1)^{i+1}]e_{-1}^{p^m-1-i}
                                         \otimes 1_\lambda\\
                                   & = & (\lambda+1)e_{-1}^{p^m-1-i}\otimes 1_\lambda\,.
\end{eqnarray*}
Since, by assumption $\lambda\ne p-1$, we have $\{e_{-1}^n\otimes 
1_\lambda\mid 1\le n\le p^m-1\}\subseteq M$. Finally, $e_1e_{-1}\otimes 1_\lambda=-\lambda
e_{-1}^0\otimes 1_\lambda$ and because $\lambda\ne 0$, one also obtains that $e_{-1}^0\otimes 
1_\lambda\in M$. Hence $\{e_{-1}^n\otimes 1_\lambda\mid 0\le n\le p^m-1\}\subseteq M$ and thus
$M=V(\lambda)$.\quad $\Box$
\bigskip

\noindent It is a consequence of Lemma \ref{ind} that most of the simple modules of the 
Zassenhaus algebras with non-vanishing cohomology are generalized baby Verma modules 
(cf.~\cite[Satz 5.7]{Leg}).

\begin{lem}\label{indirr}
If $\lambda\not\in\{0,p-1\}$, then $S(\lambda)\cong V(\lambda)$ as $W(m)$-modules.
\end{lem}

\noindent {\it Proof.\/} Set $W:=W(m)$. It follows from the above computations using the 
Cartan-Weyl identity that the linear transformation $\Phi_\lambda:V(\lambda)\to S(\lambda)$ 
given by $e_{-1}^n\otimes 1_\lambda\mapsto e_{-1}^n\cdot v$ for every $0\le n\le p^m-1$ is 
a $W$-module epimorphism. Since $\Phi_\lambda\ne 0$, and by Lemma \ref{ind} $V(\lambda)$ is 
simple, it follows that $\mathrm{Ker}(\Phi_\lambda)=0$ whence $\Phi_\lambda$ is a $W$-module 
isomorphism.\quad $\Box$
\bigskip

Now it remains to consider the two cases $\lambda=0$ and $\lambda=p-1$. Note that, by 
Proposition \ref{adjnat}(3) and (4), $V(p-1)\cong A(m)$ and $V(0)\cong A(m)^*$, respectively. 
We begin by showing that the simple module $S(0)$ is trivial (cf.~\cite[Satz 5.8]{Leg}).

\begin{lem}\label{irr0}
$S(0)$ is isomorphic to the one-dimensional trivial $W(m)$-module.
\end{lem}

\noindent {\it Proof.\/} Set $W:=W(m)$. Suppose that $e_{-1}\cdot v\ne 0$. Because $\lambda=0$,
$$e_i\cdot(e_{-1}^n\cdot v)=(-1)^{i+1}{n\choose i+1}e_{-1}^{n-i}\cdot v$$ for every $0\le
i\le p^m-2$. But ${n\choose i+1}=0$ unless $n-i>0$, so $\sum_{n=1}^{p^m-1}\mathbb{F}e_{-1}^n
\cdot v$ is a non-zero proper $W$-submodule of $S(0)$. This is a contradiction since by 
definition $S(0)$ is simple. Hence $S(0)=\mathbb{F}v$ with $e_i\cdot v=0$ every $-1\le i\le 
p^m-2$.\quad $\Box$
\bigskip

\noindent By arguments similar to those used to establish Lemma \ref{irr0}, one can prove 
that $V(0)$ has a $(p^m-1)$-dimensional submodule. The proof of Lemma \ref{ind} then shows 
that this submodule is simple (cf.~\cite[Lemma III.3.2b) for $m=1$]{Fe1}). In particular, 
$V(0)$ is indecomposable.

\begin{lem}\label{ind0}
The subspace $\bigoplus_{n=1}^{p^m-1}\mathbb{F}(e_{-1}^n\otimes 1_0)$ is a simple $W(m)$-submodule
of $V(0)$ and the factor module $V(0)/\bigoplus_{n=1}^{p^m-1}\mathbb{F}(e_{-1}^n\otimes 1_0)$ 
is isomorphic to the one-dimensional trivial $W(m)$-module.\quad $\Box$
\end{lem}

Now we consider the case $\lambda=p-1$. The following observation is crucial to understand 
the structure of $V(p-1)$ and $S(p-1)$ (cf.~\cite[Lemma 5.9]{Leg}).

\begin{lem}\label{irr-10}
If $\lambda=p-1$, then $e_{-1}^{p^m-1}\cdot v=0$, i.e., $S(p-1)=\sum_{n=0}^{p^m-2}\mathbb{F}
e_{-1}^n\cdot v$.
\end{lem}

\noindent {\it Proof.\/} Set $W:=W(m)$. Suppose that $e_{-1}^{p^m-1}\cdot v\ne 0$. Because 
$\lambda=p-1$, $$e_i\cdot(e_{-1}^n\cdot v)=(-1)^{i+1}{n+1\choose i+1}e_{-1}^{n-i}\cdot v$$ 
for every $0\le i\le p^m-2$. In particular, for $n=p^m-1$ one has that $$e_i\cdot
(e_{-1}^{p^m-1}\cdot v)=(-1)^{i+1}{p^m\choose i+1}e_{-1}^{p^m-1-i}\cdot v=0$$ for every 
$0\le i\le p^m-2$ and additionally $e_{-1}\cdot(e_{-1}^{p^m-1}\cdot v)=e_{-1}^{p^m}\cdot 
v=0$. Hence $\mathbb{F}(e_{-1}^{p^m-1}\cdot v)$ is a non-zero proper $W$-submodule of 
$S(p-1)$, a contradiction since by definition $S(p-1)$ is simple.\quad $\Box$
\bigskip

\noindent {\it Remark.\/} It follows from Lemma \ref{irr-1} below that $\dim_\mathbb{F}
S(p-1)=p^m-1$ and therefore $$S(p-1)=\bigoplus_{n=0}^{p^m-2}\mathbb{F}e_{-1}^n\cdot v\,.$$ 
\smallskip

\noindent The generalized baby Verma module $V(p-1)$ has a one-dimensional maximal submodule 
(cf.~\cite[Lemma III.3.2c)]{Fe1} for $m=1$ and \cite[Lemma 5.10]{Leg} for the general case). 
In particular, $V(p-1)$ is indecomposable.

\begin{lem}\label{ind-1}
$F:=\mathbb{F}(e_{-1}^{p^m-1}\otimes 1_{p-1})$ is a maximal $W(m)$-submodule of $V(p-1)$ 
with $W(m)F=0$. 
\end{lem}

\noindent {\it Proof.\/} Set $W:=W(m)$. As in the proof of Lemma \ref{irr-10}, one can 
show that $\mathbb{F}(e_{-1}^{p^m-1}\otimes 1_{p-1})$ is a trivial $W$-submodule of 
$V(p-1)$. Suppose that $M$ is a $W$-submodule of $V(p-1)$ that properly contains $F$ and 
let $\sum_{n=s}^{p^m-1}\alpha_ne_{-1}^n\otimes 1_{p-1}\in M$ with $\alpha_s\ne 0$ and $0
\le s\le p^m-2$. Then $$e_{-1}^{p^m-2-s}\cdot\left(\sum_{n=s}^{p^m-1}\alpha_ne_{-1}^n
\otimes 1_{p-1}\right)=\alpha_se_{-1}^{p^m-2}\otimes 1_{p-1}+\alpha_{s+1}e_{-1}^{p^m-1}
\otimes 1_{p-1}\,.$$ Since $F\subseteq M$ and $\alpha_s\ne 0$, it follows that $e_{-1}^{p^m-2}
\otimes 1_{p-1}\in M$. But $$e_{p^m-2}\cdot\left(e_{-1}^{p^m-2}\otimes 1_{p-1}\right)=
(-1)^{p^m-1}{p^m-1\choose p^m-1}e_{-1}^0\otimes 1_{p-1}=e_{-1}^0\otimes 1_{p-1}\,,$$ 
i.e., $e_{-1}^0\otimes 1_{p-1}\in M$; therefore $$e_{-1}^n\otimes 1_{p-1}=e_{-1}^n
\left(e_{-1}^0\otimes 1_{p-1}\right)\in M$$ for every $0\le n\le p^m-1$. Hence $M=V(p-1)$.
\quad $\Box$
\bigskip

\noindent Finally, we prove that $S(p-1)$ is a factor module of the generalized baby 
Verma module $ V(p-1)$ (cf.~\cite[Lemma III.3.2c)]{Fe1} for $m=1$ and \cite[Satz 5.11]{Leg} 
for the general case).

\begin{lem}\label{irr-1}
$S(p-1)\cong V(p-1)/F$ as $W(m)$-modules.
\end{lem}

\noindent {\it Proof.\/} Set $W:=W(m)$. It follows from the above computations, using 
the Cartan-Weyl identity, that the linear transformation $\Phi_{p-1}:V(p-1)\to S(p-1)$ 
given by $e_{-1}^n\otimes 1_{p-1}\mapsto e_{-1}^n\cdot v$, for every $0\le n\le p^m-1$, 
is a non-zero $W$-module epimorphism. By virtue of Lemma \ref{irr-10}, $$\Phi_{p-1}
(e_{-1}^{p^m-1}\otimes 1_{p-1})=e_{-1}^{p^m-1}\cdot v=0\,,$$ i.e., $F\subseteq\mathrm{Ker}
(\Phi_{p-1})$. Hence $\Phi_{p-1}$ induces a non-zero $W$-module epimorphism 
$\overline{\Phi}_{p-1}$ from $V(p-1)/F$ onto $S(p-1)$. According to Lemma \ref{ind-1}, 
$V(p-1)/F$ is simple and thus $\overline{\Phi}_{p-1}$ is a $W$-module isomorphism.
\quad $\Box$
\bigskip

As a consequence of Lemma \ref{ind-1} and Lemma \ref{irr-1}, we obtain the short exact 
sequence $$0\longrightarrow S(0)\longrightarrow V(p-1)\longrightarrow S(p-1)\longrightarrow 
0$$ of $W(m)$-modules. Similarly, we have the short exact sequence $$0\longrightarrow 
S(p-1)\longrightarrow V(0)\longrightarrow S(0)\longrightarrow 0$$ of $W(m)$-modules. By 
virtue of Theorem \ref{frobrec}(2), the second short exact sequence is the dual of the 
first short exact sequence. But the second short exact sequence can also be obtained 
directly by showing that the $(p^m-1)$-dimensional $W(m)$-submodule $\bigoplus_{n=1}^{p^m-1}
\mathbb{F}(e_{-1}^n\otimes 1_0)$ of $V(0)$ is isomorphic to $S(p-1)$ (see the remark 
after Lemma \ref{irr-10}).

Now we apply Theorem \ref{1cohzas} in order to determine the $1$-cohomology of the
Zassenhaus algebras with coefficients in simple modules.

\begin{thm}\label{1cohzasirr}
Let $\mathbb{F}$ be a field of prime characteristic $p>2$, let $m$ be a positive integer, 
and let $\lambda\in\mathbb{F}_p$. Then the following statements hold:
\begin{enumerate}
\item[{\rm (1)}] If $p=3$, then 
                 \begin{eqnarray*}
                 \dim_\mathbb{F}H^1(W(m),S(\lambda))=
                 \left\{
                 \begin{array}{cl}
                 0 & \mbox{if }\lambda=0\\ 
                 m-1 & \mbox{if }\lambda=1\,.\\
                 2 & \mbox{if }\lambda=2
                 \end{array}
                 \right.
                 \end{eqnarray*}
\item[{\rm (2)}] If $p>3$, then
                 \begin{eqnarray*}
                 \dim_\mathbb{F}H^1(W(m),S(\lambda))=
                 \left\{
                 \begin{array}{cl}
                 1 & \mbox{if }\lambda = 1\\ 
                 m-1 & \mbox{if }\lambda=p-2\\
                 2 & \mbox{if }\lambda=p-1\\
                 0 & \mbox{if }\lambda\ne 1,p-2,p-1\,.
                 \end{array}
                 \right.
                 \end{eqnarray*}
\end{enumerate}
\end{thm}

\noindent {\it Proof.\/} Set $W:=W(m)$ and $\mathfrak{u}:=\mathfrak{u}(m)$. Then the 
assertions for $\lambda\ne 0,p-1$ are immediate consequences of Lemma \ref{indirr} and 
Theorem \ref{1cohzas}. According to \cite[Theorem 4.2.4(1)]{SF}, $W$ is simple and 
therefore perfect which yields the assertions for $\lambda=0$. Since $H^1(W,S(0))
=0$, we obtain the exact sequence $$(++)\qquad V(0)^W\to S(0)^W\longrightarrow H^1
(W,S(p-1))\to H^1(W,V(0))\to 0$$ by applying the long exact cohomology sequence to 
the short exact sequence $$0\longrightarrow S(p-1)\longrightarrow V(0)\longrightarrow 
S(0)\longrightarrow 0\,.$$ Furthermore, it follows from Proposition \ref{invzas} that 
$V(0)^W=0$. Hence we obtain from $(++)$ that $$H^1(W,S(p-1))\cong S(0)^W\oplus H^1
(W,V(0))$$ and therefore we conclude from Lemma \ref{irr0} and Theorem \ref{1cohzas} 
that $$\dim_\mathbb{F}H^1(W,S(p-1))=\dim_\mathbb{F}S(0)^W+\dim_\mathbb{F}H^1(W,V(0))
=2$$ which finishes the proof of the theorem.\quad $\Box$
\bigskip

\noindent {\it Remark.\/} If the ground field $\mathbb{F}$ is algebraically closed, then 
Theorem \ref{1cohzasirr} completely describes the $1$-cohomology of Zassenhaus algebras 
with coefficients in simple modules.
\bigskip

In view of the discussion before Lemma \ref{triv}, every simple $W(m)$-module is also 
simple as a $\mathfrak{W}(m)$-module. In particular, $S(\lambda)$ is a simple $\mathfrak{W}
(m)$-module for every $\lambda\in\mathbb{F}_p$. The following result is an immediate 
consequence of Corollary \ref{lowcohpenv}(2) and Theorem \ref{1cohzasirr}:

\begin{thm}\label{1cohreszasirr}
Let $\mathbb{F}$ be a field of prime characteristic $p>2$, let $m$ be a positive integer, 
and let $\lambda\in\mathbb{F}_p$. Then the following statements hold:
\begin{enumerate}
\item[{\rm (1)}] If $p=3$, then 
                 \begin{eqnarray*}
                 \dim_\mathbb{F}H^1(\mathfrak{W}(m),S(\lambda))=
                 \left\{
                 \begin{array}{cl}
                 m-1 & \mbox{if }\lambda=0\mbox{ or }\lambda=1\\
                 2 & \mbox{if }\lambda=2\,.
                 \end{array}
                 \right.
                 \end{eqnarray*}
\item[{\rm (2)}] If $p>3$, then
                 \begin{eqnarray*}
                 \dim_\mathbb{F}H^1(\mathfrak{W}(m),S(\lambda))=
                 \left\{
                 \begin{array}{cl}
                 m-1 & \mbox{if }\lambda=0\mbox{ or }\lambda=p-2\\
                 1 & \mbox{if }\lambda=1\\
                 2 & \mbox{if }\lambda=p-1\\
                 0 & \mbox{if }\lambda\ne 0,1,p-2,p-1\,.\quad\Box
                 \end{array}
                 \right.
                 \end{eqnarray*}
\end{enumerate}
\end{thm}

\noindent {\it Remark.\/} If $S$ is a simple $\mathfrak{W}(m)$-module with $H^\bullet
(\mathfrak{W}(m),S)=0$, then it follows from Corollary \ref{vancohirred} that $S$ is 
a restricted $\mathfrak{W}(m)$-module. One can show that $S$ is isomorphic to $S(\lambda)$ 
for some $\lambda\in\mathbb{F}_p$ by employing arguments similar to those used above in 
the case of $W(m)$ and by replacing the generalized baby Verma module $V(\lambda)$ by 
the restricted baby Verma module $u(\mathfrak{W}(m))\otimes_{u(\mathfrak{b}(m))}F_\lambda$. 
As a consequence, Theorem \ref{1cohreszasirr} completely describes the $1$-cohomology 
of $\mathfrak{W}(m)$ with coefficients in simple modules if the ground field $\mathbb{F}$ 
is algebraically closed.
\bigskip

\noindent Now Theorem \ref{1cohreszasirr} can be used to determine the dimensions of 
the restricted $1$-cohomology of the minimal $p$-envelopes of Zassenhaus algebras with 
coefficients in simple restricted modules (cf.~\cite[Satz III.3.3]{Fe1} for $m=1$ and 
$p>3$).

\begin{cor}\label{1rescohreszasirr}
Let $\mathbb{F}$ be a field of prime characteristic $p>2$, let $m$ be a positive integer, 
and let $\lambda\in\mathbb{F}_p$. Then the following statements hold:
\begin{enumerate}
\item[{\rm (1)}] If $p=3$, then
                 \begin{eqnarray*}
                 \dim_\mathbb{F}H_*^1(\mathfrak{W}(m),S(\lambda))=
                 \left\{
                 \begin{array}{cl}
                 0 & \mbox{if }\lambda=0\\
                 m-1 & \mbox{if }\lambda=1\\
                 2 & \mbox{if }\lambda=2\,.
                 \end{array}
                 \right.
                 \end{eqnarray*}
\item[{\rm (2)}] If $p>3$, then
                 \begin{eqnarray*}
                 \dim_\mathbb{F}H_*^1(\mathfrak{W}(m),S(\lambda))=
                 \left\{
                 \begin{array}{cl}
                 1 & \mbox{if }\lambda=1\\ 
                 m-1 & \mbox{if }\lambda=p-2\\
                 2 & \mbox{if }\lambda=p-1\\
                 0 & \mbox{if }\lambda\neq 1,p-1,p-2\,.
                 \end{array}
                 \right.
                 \end{eqnarray*}
\end{enumerate}
\end{cor} 

\noindent {\it Proof.\/} The vanishing of $H_*^1(\mathfrak{W}(m),S(0))$ follows from 
\cite[Proposition 2.7]{Fe2}, the perfectness of $W(m)$, and the definition of the 
$p$-mapping of $\mathfrak{W}(m)$ (see Section 1). The remaining statements can be 
obtained from Theorem \ref{1cohreszasirr} and Proposition \ref{invzas} in conjunction 
with the six-term exact sequence of Hochschild \cite[p.~575]{Ho} or the Appendix of 
\cite{Fe3}.\quad $\Box$
\bigskip

We conclude the paper by determining the central extensions of the Zassenhaus algebras and 
their minimal $p$-envelopes. For $p>3$ the next result is due to Block \cite[Theorem 
5.1]{Blo} and the corresponding non-trivial central extensions of the Zassenhaus algebras 
are the {\it modular Virasoro algebras\/} (cf.~\cite{Dzh3}). The case $p=3$ is due to 
Dzhumadil'daev \cite[Theorem 2]{Dzh2} (cf.~also \cite[Theorem 3.2]{Fa2} and \cite[Theorem 
3.1 and Corollary 3.1]{Chiu}).

\begin{thm}\label{2cohzas}
Let $\mathbb{F}$ be a field of prime characteristic $p$ and let $m$ be a positive integer. 
Then the following statements hold:
\begin{enumerate}
\item[{\rm (1)}] If $p=3$, then $\dim_\mathbb{F}H^2(W(m),\mathbb{F})=m-1$.
\item[{\rm (2)}] If $p>3$, then $\dim_\mathbb{F}H^2(W(m),\mathbb{F})=1$.
\end{enumerate}
\end{thm}

\noindent {\it Proof.\/} Set $W:=W(m)$ and $A:=A(m)$.

(1): If $p=3$, then it follows from Proposition \ref{adjnat}(2) and Theorem \ref{1cohzas}(2)
that $$\dim_\mathbb{F}H^1(W,W^*)=\dim_\mathbb{F}H^1(W,V(1))=m-1\,.$$ According to the remark
after Proposition \ref{adjnat}, $W^*\cong W$ as $W$-modules. For any $0\le a\le 3^m-1$ let 
$\gamma_a:A\to\mathbb{F}$ be the linear transformation defined by $$f=\sum_{a=0}^{3^m-1}
\gamma_a(f)x^{(a)}\qquad\forall~f\in A\,.$$ Then the canonical pairing of $W^*$ and $W$ 
induces via the above isomorphism a non-degenerate invariant bilinear form on $W$ given by 
$(f\partial,g\partial):=\gamma_{2^m-1}(fg)$ (cf.~the proof of \cite[Theorem 4.6.3]{SF}). It 
follows from the remark after Theorem \ref{1cohzas} that the $m-1$ linear independent outer 
derivations of $W$ are induced by $\partial^{3^r}$ for $1\le r\le m-1$. Then 
\begin{eqnarray*}
([\partial^{3^r},f\partial],g\partial)+(f,[\partial^{3^r},g\partial]) 
& = & \gamma_{3^m-1}(\partial^{3^r}(f)g+f\partial^{3^r}(g))\\ 
& = & \gamma_{3^m-1}(\partial^{3^{r-1}+1}(\partial^2(f)g-\partial(f)\partial(g)+f\partial^2(g)))\\
& = & 0
\end{eqnarray*}
and in the light of the discussion before \cite[Proposition 1]{Dzh2}, we conclude that 
$\dim_\mathbb{F}H^2(W,\mathbb{F})=m-1$.

(2): If $p>3$, then it follows from Proposition \ref{adjnat}(2) and Theorem \ref{1cohzas}(3)
that $$\dim_\mathbb{F}H^1(W,W^*)=\dim_\mathbb{F}H^1(W,V(1))=1$$ and the assertion can be
obtained from \cite[Proposition 1.3(3)]{Fa1} in conjunction with \cite[Theorem 4.2.4(1)]{SF} 
and the remark after Proposition \ref{adjnat}.\quad $\Box$
\bigskip

\noindent {\it Remark.\/} In \cite[Theorem 2]{Dzh3} Dzhumadil'daev also claims to have determined 
the central extensions of the Zassenhaus algebras for $p=2$ (cf.~also \cite[Theorem 3.1 and 
Corollary 3.1]{Chiu}). But at least in the case $m=1$ this result is not correct since for 
$p=2$ the Witt algebra is the two-dimensional non-abelian Lie algebra which has {\it no\/} 
non-trivial central extensions.
\bigskip

Finally, as an immediate consequence of Corollary \ref{lowcohpenv}(3) and Theorem \ref{2cohzas} 
we find the central extensions of the minimal $p$-envelopes of the Zassenhaus algebras.

\begin{thm}\label{2cohreszas}
Let $\mathbb{F}$ be a field of prime characteristic $p>2$ and let $m$ be a positive integer. 
Then the following statements hold:
\begin{enumerate}
\item[{\rm (1)}] If $p=3$, then $\dim_\mathbb{F}H^2(\mathfrak{W}(m),\mathbb{F})=\frac{(m-1)m}{2}$.
\item[{\rm (2)}] If $p>3$, then $\dim_\mathbb{F}H^2(\mathfrak{W}(m),\mathbb{F})=\frac{m^2-3m+4}{2}$.
                 \quad $\Box$
\end{enumerate}
\end{thm}

\noindent {\it Remark.\/} From Theorem \ref{2cohreszas} and \cite[Proposition 1]{Dzh2} one 
obtains lower bounds for the dimensions of the first coadjoint cohomology $H^1(\mathfrak{W}
(m),\mathfrak{W}(m)^*)$ of the minimal $p$-envelopes of the Zassenhaus algebras.
\bigskip

\noindent We leave it to the interested reader to apply Theorem \ref{2cohreszas} in conjunction 
with the six-term exact sequence of Hochschild \cite[p.~575]{Ho} to obtain upper bounds for the 
dimensions of the second restricted cohomology of $\mathfrak{W}(m)$ with coefficients in the 
one-dimensional trivial module. It would be interesting to know the exact dimensions of
$H_*^2(\mathfrak{W}(m),\mathbb{F})$.




\begin{thebibliography}{99}


\bibitem{Blo}
R. E. Block: 
On the extensions of Lie algebras,
{\it Canad. J. Math.\/} {\bf 20} (1968), 1439--1450.

\bibitem{Bou}
N. Bourbaki:
{\it Lie Groups and Lie Algebras: Chapters 1--3\/} (2$^{\rm nd}$ printing),
Springer-Verlag, Berlin/Heidelberg/New York/London/Paris/Tokyo, 1989.

\bibitem{Chiu}
Sen Chiu: 
Central extensions and $H^1(L,L^*)$ of the graded Lie algebras of Cartan type,
{\it J. Algebra\/} {\bf 149} (1992), 46--67.

\bibitem{CS}
Sen Chiu and Guangyu Shen: 
Cohomology of graded Lie algebras of Cartan type of characteristic $p$,
{\it Abh. Math. Sem. Univ. Hamburg\/} {\bf 57} (1987), 139--156.

\bibitem{Dix}
J. Dixmier:
Cohomologie des alg\`ebres de Lie nilpotentes,
{\it Acta Sci. Math. (Szeged)\/} {\bf 16} (1955), 246--250.

\bibitem{Dzh1}
A. S. Dzhumadil'daev:
On the cohomology of modular Lie algebras,
{\it Math. USSR Sb.\/} {\bf 47} (1984), 127--143. 

\bibitem{Dzh2}
A. S. Dzhumadil'daev:
Central extensions and invariant forms of Cartan type Lie algebras of positive characteristic,
{\it Funct. Anal. Appl.\/} {\bf 18} (1984), 331--332.

\bibitem{Dzh3}
A. S. Dzhumadil'daev:
Central extensions of the Zassenhaus algebra and their irreducible representations,
{\it Math. USSR Sb.\/} {\bf 54} (1986), 457--474. 

\bibitem{Dzh4}
A. S. Dzhumadil'daev:  
Cohomology of truncated coinduced representations of Lie algebras of positive characteristic,
{\it Math. USSR Sb.\/} {\bf 66} (1990), 461--473. 

\bibitem{Dzh5}
A. S. Dzhumadil'daev:  
Cohomology and nonsplit extensions of modular Lie algebras,
in: {\it Proceedings of the International Conference on Algebra Dedicated to the Memory of 
A. I. Mal'cev, Novosibirsk, 1989\/} (eds. L. A. Bokut', Yu. L. Ershov, and A. I. Kostrikin), 
Contemp. Math., vol.~{\bf 131}, Part 2, Amer. Math. Soc., Providence, RI, 1992, pp.~31--43.

\bibitem{Fa1}
R. Farnsteiner:
Central extensions and invariant forms of graded Lie algebras,
{\it Algebras Groups Geom.\/} {\bf 3} (1986), 431--455. 

\bibitem{Fa2}
R. Farnsteiner:
Dual space derivations and $H^2(L,F)$ of modular Lie algebras,
{\it Canad. J. Math.\/} {\bf 39} (1987), 1078--1106.

\bibitem{Fa3}
R. Farnsteiner:
On the cohomology of associative algebras and Lie algebras,
{\it Proc. Amer. Math. Soc.\/} {\bf 99} (1987), 415--420.

\bibitem{Fa4}
R. Farnsteiner:
On the vanishing of homology and cohomology groups of associative algebras,
{\it Trans. Amer. Math. Soc.\/} {\bf 306} (1988), 651--665.

\bibitem{Fa5}
R. Farnsteiner:
Extension functors of modular Lie algebras,
{\it Math. Ann.\/} {\bf 288} (1990), 713--730.

\bibitem{FS}
R. Farnsteiner and H. Strade:
Shapiro's lemma and its consequences in the cohomology theory of modular Lie algebras,
{\it Math. Z.\/} {\bf 206} (1991), 153--168.

\bibitem{Fe1}
J. Feldvoss:
{\it Homological Aspects in the Representation Theory of Modular Lie Algebras\/} (in German), 
Doctoral Dissertation, University of Hamburg, 1989. 

\bibitem{Fe2}
J. Feldvoss:
On the cohomology of restricted Lie algebras, 
{\it Comm. Algebra\/} {\bf 19} (1991), 2865--2906.

\bibitem{Fe3}
J. Feldvoss:
A cohomological characterization of solvable modular Lie algebras, 
in: {\it Non-Associative Algebra and Its Applications, Oviedo, 1993\/} (ed. S. Gonz\'alez), 
Mathematics and Its Applications, vol.~{\bf 303}, 
Kluwer Academic Publishers, Dordrecht/Boston/London, 1994, pp.~133--139.

\bibitem{Ho}
G. Hochschild:
Cohomology of restricted Lie algebras,
{\it Amer. J. Math.\/} {\bf 76} (1954), 555--580.

\bibitem{HS}
G. Hochschild and J.-P. Serre: 
Cohomology of Lie algebras,
{\it Ann. of Math.\/} (2) {\bf 57} (1953), 591--603.

\bibitem{Jac} 
N. Jacobson: 
{\it Lie Algebras\/}, 
Dover Publications, Inc., New York, 1979 
(unabridged and corrected republication of the original edition from 1962).

\bibitem{Kuz}
M. I. Kuznetsov: 
Truncated induced modules over transitive Lie algebras of characteristic $p$, 
{\it Math. USSR Izv.\/} {\bf 34} (1990), 575--608.

\bibitem{Leg}
K. Legler:
{\it Cohomology Groups of Modular Lie Algebras\/} (in German),
Diploma Thesis, University of Hamburg, 1990.

\bibitem{Mil}
A. A. Mil'ner:  
Irreducible representations of Zassenhaus algebras (in Russian), 
{\it Uspekhi Mat. Nauk\/} {\bf 30} (1975), 178.

\bibitem{Shu}
Bin Shu: 
Generalized restricted Lie algebras and representations of the Zassenhaus algebra,
{\it J. Algebra\/} {\bf 204} (1998), 549--572.

\bibitem{Str1}
H. Strade:
The role of $p$-envelopes in the theory of modular Lie algebras,
in: {\it Lie Algebras and Related Topics, Madison, WI, 1988\/} (eds. G. Benkart and J. M. Osborn), 
Contemp. Math., vol.~{\bf 110}, Amer. Math. Soc., Providence, RI, 1990, pp.~265--287. 

\bibitem{Str2}
H. Strade:
{\it Simple Lie Algebras over Fields of Positive Characteristic I: Structure Theory\/},
de Gruyter Expositions in Mathematics, vol.~{\bf 38},
Walter de Gruyter \& Co., Berlin/New York, 2004.

\bibitem{SF}  
H. Strade and R. Farnsteiner: 
{\it Modular Lie Algebras and Their Representations\/}, 
Monographs and Textbooks in Pure and Applied Mathematics, vol.~{\bf 116}, 
Marcel Dekker, Inc., New York/Basel, 1988.

\bibitem{Wei}
C. A. Weibel:
{\it An Introduction to Homological Algebra\/},
Cambridge Studies in Advanced Mathematics, vol.~{\bf 38},
Cambridge University Press, Cambridge, 1994.


\end{thebibliography}
\end{document}